\RequirePackage{ifpdf}
\ifpdf 
\documentclass[pdftex]{sigma}
\else
\documentclass{sigma}
\fi

\numberwithin{equation}{section}
\numberwithin{theorem}{section}

\def\fratop{\genfrac{}{}{0pt}1}
\def\satop#1#2{\fratop{\scriptstyle#1}{\scriptstyle#2}}

\let\leq\leqslant
\let\ge\geqslant
\let\geq\geqslant
\let\alb\allowbreak

\def\>{\relax\ifmmode\mskip.666667\thinmuskip\relax\else\kern.111111em\fi}
\def\<{\relax\ifmmode\mskip-.333333\thinmuskip\relax\else\kern-.0555556em\fi}
\def\id{\mathrm{id}}
\def\vsk#1>{\vskip#1\baselineskip}
\def\vvn#1>{\vadjust{\vsk#1>}}

\def\diag{\mathop{\mathrm{diag}}\nolimits}
\def\Sym{\mathop{\mathrm{Sym}}\nolimits}
\def\XXX/{{\slshape XXX}}

\newtheorem{lemma}[theorem]{Lemma}
\newtheorem{corollary}[theorem]{Corollary}
\newtheorem{proposition}[theorem]{Proposition}
{\theoremstyle{definition}
\newtheorem{example}[theorem]{Example}
\newtheorem{remark}[theorem]{Remark}
}

\newcommand{\nc}{\newcommand}
\nc{\on}{\operatorname}
\nc{\ch}{\mbox{ch}}
\nc{\Z}{{\mathbb Z}}
\nc{\C}{{\mathbb C}}
\nc{\R}{{\mathbb R}}
\nc{\pone}{{\mathbb C}{\mathbb P}^1}
\nc{\p}{\partial_u}
\nc{\F}{{\mathcal F}}
\nc{\arr}{\rightarrow}
\nc{\larr}{\longrightarrow}
\nc{\al}{\alpha}
\nc{\ri}{\rangle}
\nc{\lef}{\langle}
\nc{\W}{{\mathcal W}}
\nc{\la}{\lambda}
\nc{\ep}{\epsilon}
\nc{\eps}{\varepsilon}

\nc{\Om}{\Omega}
\nc{\su}{\widehat{{\mathfrak sl}}_2}
\nc{\sw}{{\mathfrak s}{\mathfrak l}}

\nc{\g}{{\mathfrak g}}
\nc{\h}{{\mathfrak h}}
\nc{\n}{{\mathfrak n}}
\nc{\N}{\widehat{\n}}
\nc{\G}{\widehat{\g}}
\nc{\De}{\Delta_+}
\nc{\gt}{\widetilde{\g}}
\nc{\Ga}{\Gamma}
\nc{\one}{{\mathbf 1}}
\nc{\z}{{\mathfrak Z}}
\nc{\zz}{{\mathcal Z}}
\nc{\Hh}{{\mathcal H}_\beta}
\nc{\qp}{q^{\frac{k}{2}}}
\nc{\qm}{q^{-\frac{k}{2}}}
\nc{\La}{\Lambda}
\nc{\wt}{\widetilde}
\nc{\qn}{\frac{[m]_q^2}{[2m]_q}}
\nc{\cri}{_{\on{cr}}}
\nc{\kk}{h^\vee}
\nc{\sun}{\widehat{\sw}_N}
\nc{\hh}{\widehat{\mathfrak h}}
\nc{\HH}{{\mathcal H}_{q,t}}
\nc{\ca}{\wt{{\mathcal A}}_{h,k}(\sw_2)}

\nc{\gl}{\widehat{{\mathfrak g}{\mathfrak l}}_2}
\nc{\el}{\ell}
\nc{\s}{{\mathbf s}}
\nc{\bi}{\bibitem}
\nc{\om}{\omega}
\nc{\WW}{\W_\beta}
\nc{\scr}{{\mathbf S}}
\nc{\ab}{{\mathbf a}}
\nc{\rr}{r}
\nc{\ol}{\overline}
\nc{\con}{qt^{-1} + q^{-1}t}
\nc{\den}{q^{\el-1} t^{-\el+1}+ q^{-\el+1} t^{\el-1}}
\nc{\ds}{\displaystyle}
\nc{\B}{B}
\nc{\A}{{\mathbb A}}
\nc{\GG}{{\mathcal G}}
\nc{\UU}{{\mathcal U}}
\nc{\MM}{{\mathcal M}}
\nc{\CC}{{\mathcal C}}
\nc{\GL}{{}^L G}
\nc{\dzz}{\frac{dz}{z}}
\nc{\Res}{\on{Res}}
\nc{\rep}{{\mathcal R}ep \;}
\nc{\uqg}{U_q \G}
\nc{\uqgg}{U_q \g}
\nc{\Fq}{{\mathbb F}_q}

\nc{\stimes}{\ltimes}
\nc{\K}{\hat{\mathcal K}}
\nc{\Ql}{\ol{\mathbb Q}_\ell}

\nc{\ga}{\gamma}
\nc{\PL}{{}^L P}
\nc{\E}{\mc E}
\nc{\mc}{\mathcal}
\nc{\mbf}{\mathbf}
\nc{\bb}{{\mathfrak b}}
\nc{\OO}{{\mc O}}
\nc{\Po}{{\mc P}}
\nc{\V}{{\mc V}}
\nc{\yy}{{\mc Y}}
\nc{\M}{\mathcal M}
\nc{\Coh}{{{\mathcal C}oh}}
\nc{\Cohn}{\Coh_n}
\nc{\f}{{\mathcal F}}
\nc{\si}{_E}
\nc{\Gaf}{{\mathbb G}_{a,\Fq}}
\nc{\KK}{{\mathfrak k}}

\nc{\PO}{{\mathbb P^1}}
\nc{\PR}{{\mathbb P^r}}

\nc{\Wr}{{ {\rm Wr}}}

\nc{\bs}{\boldsymbol}
\nc{\Ref}[1]{{\rm(\ref{#1})}}

\nc{\U}{\mathcal{U}}

\nc{\sing}{{\rm Sing}\,}

\nc{\Ll}{\La - \al(\bs l)}

\nc{\nL}{L_{\om_r}^{\otimes n}[\mu]}
\nc{\nnL}{L_{\om_r}^{\otimes n}}

\nc{\snL}{ \sing L_{\om_r}^{\otimes n}[\mu]}
\nc{\btz}{ \om(\bs t; \bs z)}
\nc{\SL}{ \mathfrak{sl}}
\nc{\GR}{ {G(r+1,d)}}

\def\glN{\mathfrak{gl}_N}

\def\Dg{\mathcal D}
\let\si\sigma

\nc{\Te}{\mathfrak{T}}

\nc{\End}{ {\rm End}}
\nc{\Wrd}{{\rm Wr}^{\rm(d)}}

\nc{\Y}{Y(\glN)}
\nc{\MLz}{M_{\bs \La}(\bs z)}

\def\Ys{Y(\mathfrak{sl}_N)}
\nc{\D}{\mathfrak {D}}

\nc{\DM}{\D_{Q, M(\bs z)}(u,\tau)}

\nc{\Dt}{\D_{\tilde {\bs t}}(u,\tau)}
\nc{\Se}{\mathfrak S}
\nc{\Xe}{\mathcal X}
\nc{\Tee}{\mathcal S}
\nc{\glx}{\glN[x]}
\nc{\dd}{\mathfrac d{du}}

\begin{document}

\allowdisplaybreaks

\renewcommand{\PaperNumber}{060}

\FirstPageHeading

\renewcommand{\thefootnote}{$\star$}

\ShortArticleName{Generating Operator of \XXX/ or Gaudin Transfer
Matrices Has Quasi-Exponential Kernel}

\ArticleName{Generating Operator of \XXX/ or Gaudin Transfer\\
Matrices Has Quasi-Exponential Kernel\footnote{This paper is a
contribution to the Vadim Kuznetsov Memorial Issue `Integrable
Systems and Related Topics'. The full collection is available at
\href{http://www.emis.de/journals/SIGMA/kuznetsov.html}{http://www.emis.de/journals/SIGMA/kuznetsov.html}}}

\Author{Evgeny MUKHIN~$^\dag$, Vitaly TARASOV~$^{\dag\ddag}$
and Alexander~VARCHENKO~$^\S$}

\AuthorNameForHeading{E.~Mukhin, V.~Tarasov and A.~Varchenko}

\Address{$^\dag$~Department of Mathematical Sciences,
Indiana University\,--\,Purdue University Indianapolis,\\
$\phantom{^\dag}$~402 North Blackford St, Indianapolis, IN 46202-3216, USA}
\EmailD{\href{mailto:mukhin@math.iupui.edu}{mukhin@math.iupui.edu}, \href{mailto:vtarasov@math.iupui.edu}{vtarasov@math.iupui.edu}}

\Address{$^\ddag$~St.\,Petersburg Branch of Steklov Mathematical Institute,\\
$\phantom{^\ddag}$~Fontanka 27, St.\,Petersburg, 191023, Russia}
\EmailD{\href{mailto:vt@pdmi.ras.ru}{vt@pdmi.ras.ru}}

\Address{$^\S$~Department of Mathematics, University of North Carolina
at Chapel Hill,\\
$\phantom{^\S}$~Chapel Hill, NC 27599-3250, USA}
\EmailD{\href{mailto:anv@email.unc.edu}{anv@email.unc.edu}}

\ArticleDates{Published online April 25, 2007; Misprints corrected: March 2013}

\Abstract{Let $M$ be the tensor product of f\/inite-dimensional polynomial evaluation
$\Y$-modules. Consider the universal dif\/ference operator
$\D = \sum\limits_{k=0}^N (-1)^k \Te_k(u) e^{-k\partial _u }$ whose coef\-f\/icients
$\Te_k(u): M \to M$ are the \XXX/ transfer matrices associated with $M$.
We show that the dif\/ference equation $\D f = 0$ for an $M$-valued
function $f$ has a basis of solutions consisting of quasi-exponentials.
We prove the same for the universal dif\/ferential operator
$D = \sum\limits_{k=0}^N (-1)^k\Tee_k(u) \partial_u^{N-k}$ whose coef\/f\/icients
$\Tee_k(u) : \M \to \M$ are the Gaudin transfer matrices associated with the
tensor product $\M$ of f\/inite-dimensional polynomial evaluation $\glx$-modules.}

\Keywords{Gaudin model; \XXX/ model; universal dif\/ferential operator}

\Classification{34M35; 82B23; 17B67}

\begin{flushright}
\it To the memory of Vadim Kuznetsov
\end{flushright}

\section{Introduction}

In a quantum integrable model one constructs a collection of one-parameter
families of commu\-ting linear operators (called transfer matrices) acting on
a f\/inite-dimensional vector space. In this paper we consider the vector-valued
dif\/ferential or dif\/ference linear operator whose coef\/f\/icients are the transfer
matrices. This operator is called the universal dif\/ferential or dif\/ference
operator. For the \XXX/ and Gaudin type models we show that the kernel of the
universal operator is generated by quasi-exponentials or sometimes just
polynomials. This statement establishes a~relationship between these quantum
integrable models and that part of algebraic geometry, which studies the
f\/inite-dimensional spaces of quasi-exponentials or polynomials, in particular
with Schubert calculus.

We plan to develop this relationship in subsequent papers.
An example of an application of this relationship see in \cite{MTV1}.

Consider the complex Lie algebra $\glN$ with standard generators $e_{ab}$,
$a,b=1,\dots,N$. Let $\bs \La = (\La^{(1)},\dots,\La^{(n)})$ be a collection
of integral dominant $\glN$-weights, where
$\La^{(i)}=(\La_1^{(i)},\dots, \La_N^{(i)})$ for $i=1,\dots, n$
and $\La^{(i)}_N \in \Z_{\geq 0}$. Let
$M \ =\ M_{\La^{(1)}} \otimes \dots \otimes M_{\La^{(n)}}$
be the tensor product of the corresponding f\/inite-dimensional
highest weight $\glN$-modules.
Let $K = (K_{ab})$ be an $N\times N$ matrix with complex entries.
Let $z_1,\dots,z_n$ be distinct complex numbers. For $a,b=1,\dots, N$, set
\[
X_{ab}(u,\p) = \delta_{ab} \,\p
- K_{ab} - \sum_{j=1}^n\ \frac{e^{(j)}_{ba}}{u-z_j} ,
\]
where $\p=d/du$. Following \cite{T}, introduce the dif\/ferential operator
\[
\Dg_{K,M}(u,\p)=
\sum_{\si\in S_N} (-1)^\si\,
X_{1\,\si_1}(u,\p)\,X_{2\,\si_2}(u,\p)\,\cdots
\,X_{N\si_N}(u,\p) .
\]
The operator $\Dg_{K,M}(u,\p)$ acts on $M$-valued functions in $u$.
The operator is called {\it the universal
differential operator} associated with $M$, $K$ and $z_1,\dots,z_n$.

Introduce the coef\/f\/icients $\Tee_{0,K,M}(u),\dots, \Tee_{N,K,M}(u)$:
\[
\Dg_{K,M} (u,\p) = \sum_{k=0}^N\,(-1)^k\,\Tee_{k,K,M}(u)\,\p^{N-k} .
\]
The coef\/f\/icients are called {\it the transfer matrices of the Gaudin model}
associated with $M$, $K$ and $z_1,\dots,z_n$.
The transfer matrices form a commutative family:
$[\Tee_{k,K,M}(u),\Tee_{l,K,M}(v)] = 0$ for all $k$, $l$, $u$, $v$.

A quasi-exponential is a f\/inite sum of functions of the form $e^{\la u}p(u)$,
where $p(u)$ is a polynomial.

\begin{theorem}
The kernel of the universal differential operator $\Dg_{K,M} (u,\p)$ is
generated by $M$-valued quasi-exponentials. Possible exponents $\la$ appearing
in the kernel are eigenvalues of the matrix $K$. In particular, if $K=0$, then
the kernel is generated by $M$-valued polynomials.
\end{theorem}

For $K=0$, the theorem was conjectured in \cite{CT} and proved in \cite{MTV1}.

\medskip

We describe possible degrees of polynomials $p(u)$ appearing in the kernel
further in the paper.

The universal dif\/ferential operator has singular points at $z_1,\dots,z_n$.
We describe behavior of elements of the kernel at these points.

The tensor product $M$ may be naturally regarded as the tensor product
of polynomial evalua\-tion modules over the current algebra $\glx$. Then the
operators $X_{ab}(u)$, $\Dg_{K,M}(u,\p)$ may be naturally def\/ined in terms
of the $\glx$-action on $M$.

Similarly, the tensor product $M$ may be naturally regarded as the tensor
product of polynomial evaluation modules over the Yangian $\Y$. Then one may
def\/ine a linear \hbox{$N$-}th order dif\/ference operator acting on $M$-valued
functions in $u$. That operator is called the universal dif\/ference operator.
Its coef\/f\/icients commute and are called the transfer matrices of the associa\-ted
\XXX/ model.

We prove that the kernel of the universal dif\/ference operator is generated by
quasi-expo\-nentials. We describe the quasi-exponentials entering the kernel
and their behavior at singular points of the universal dif\/ference operator.

The paper has the following structure.
In Section \ref{sec quasi-exponentials}
we make general remarks on quasi-exponentials.

In Section \ref{Sec Generating} we collect basic facts about the Yangian $\Y$,
the fundamental dif\/ference operator and the \XXX/ transfer matrices.
In Section \ref{The first main} we formulate
Theorem \ref{The first main thm} which states that the kernel of the universal
dif\/ference operator is generated by quasi-exponentials.
Theorem~\ref{The first main thm} is our f\/irst main result for the \XXX/ type
models.

In Section~\ref{Sec Cont} we prove a continuity principle for dif\/ference
operators with quasi-exponen\-tial kernel. Under certain conditions we show
that if a family of dif\/ference operators has a~limiting dif\/ference operator,
and if the kernel of each operator in the family is generated by
quasi-exponentials, then the kernel of the limiting dif\/ference operator
is generated by quasi-exponentials too.

Section \ref{Bethe ansatz} is devoted to the Bethe ansatz method for the \XXX/
type models. Using the Bethe ansatz method we prove the special case of
Theorem \ref{The first main thm} in which $M$ is the tensor product of vector
representations of $\glN$. Then the functoriality properties of the fundamental
dif\/ference operator and our continuity principle allow us to deduce the general
case of Theorem \ref{The first main thm} from the special one.

In Section \ref{Sec Comparison theorem} we give a formula comparing the kernels
of universal dif\/ference operators associated respectively with the tensor
products $M_{\La^{(1)}} \otimes \dots \otimes M_{\La^{(n)}}$
and $M_{\La^{(1)}(a_1)} \otimes \dots \otimes M_{\La^{(n)}(a_n)}$,
where $a_1,\dots,a_n$ are non-negative integers and
$\La^{(i)}(a_i) = \big(\La^{(i)}_1+a_i,\dots,\La^{(i)}_N+a_i\big)$.

In Section \ref{The kernel of D} we describe the quasi-exponentials entering
the kernel of the universal dif\/ference operator and the behavior of functions
of the kernel at singular points of the universal dif\/ference operator.
Theorems \ref{thm second main, part 1}, \ref{sec local result} and
\ref{thm on space for Q=1} form our second main result for the \XXX/ type
models.

In Sections \ref{Sec Gen Gaudin}--\ref{The kernel of D Gaudin} we develop
an analogous theory for the universal dif\/ferential operator of the Gaudin type
models.

Section \ref{Sec Gen Gaudin} contains basic facts about the current algebra
$\glx$, the universal dif\/ferential operator and the Gaudin transfer matrices.
We formulate Theorem \ref{The first main thm Gaudin}, which states that
the kernel of the universal dif\/ferential operator is generated
by quasi-exponentials. This theorem is our f\/irst main result
for the Gaudin type models.

Section \ref{Sec Contin Gaudin} contains the continuity principle
for the dif\/ferential operators with quasi-exponential kernel.

Section \ref{Bethe ansatz Gaudin} is devoted to the Bethe ansatz method
for the Gaudin type model.

In Section \ref{Comparison theorem in the Gaudin case} we compare kernels of
the universal dif\/ferential operators associated respectively with tensor
products $M_{\La^{(1)}} \otimes \dots \otimes M_{\La^{(n)}}$ and
$M_{\La^{(1)}(a_1)} \otimes \dots \otimes M_{\La^{(n)}(a_n)}$.

In Section \ref{The kernel of D Gaudin} we describe the quasi-exponentials
entering the kernel of the universal dif\/fer\-ential operator and the behavior
of functions of the kernel at singular points of the universal dif\/ferential
operator. Theorems \ref{thm second main, part 1 Gaudin},
\ref{sec local result Gaudin} and \ref{thm on space for K=0} form our second
main result for the Gaudin type models.

\section{Spaces of quasi-exponentials}
\label{sec quasi-exponentials}

\subsection{Quasi-exponentials}
{\bf 2.1.1.}
Def\/ine the operator $\tau$ acting on functions of $u$ as
$(\tau f)(u)=f(u+1)$.
A function $f(u)$ will be called {\it one-periodic} if $f(u+1) = f(u)$.
Meromorphic one-periodic functions form a f\/ield with respect to addition
and multiplication.


\medskip
\noindent
{\bf 2.1.2.}
Let $Q$ be a nonzero complex number with f\/ixed argument.
Set $Q^u = e^{u\ln Q}$. We have $\tau Q^u = Q^u Q$.

Let $p\in\C[u]$ be a polynomial.
The function $Q^u p $ will be called {\it a (scalar) elementary
quasi-exponential in $u$}.
A f\/inite sum of elementary quasi-exponentials will be called
{\it a (scalar) quasi-exponential}.

Let $V$ be a complex vector space of f\/inite dimension $d$.
{\it A $V$-valued quasi-exponential} is a $V$-valued function
of the form $\sum_{a} f_a(u) v_{a}$, where $f_a(u)$ are scalar
quasi-exponentials, $v_{a}\in V$, and the sum is f\/inite.

We say that a quasi-exponential $\sum_{ab} Q_a^{u} u^b v_{ab}$
is of degree less than $k$ if $v_{ab} = 0$ for all $b\geq k$.


\medskip

\noindent
{\bf 2.1.3.} For given $\End\,(V)$-valued rational functions
$A_0(u),\dots, A_N(u)$ consider the dif\/ference operator
\begin{gather}
\label{diff oper}
\D = \sum_{k=0}^N A_k(u)\, \tau^{-k}
\end{gather}
acting on $V$-valued functions in $u$.

We say that the kernel of $\D$ {\it is generated by quasi-exponentials}
if there exist $N d$ quasi-exponential functions with values in $V$ such that
each of these function belongs to the kernel of~$\D$ and these functions
generate an $Nd$-dimensional vector space over the f\/ield of one-periodic
meromorphic functions.

The following simple observation is useful.

\begin{lemma}
Assume that a quasi-exponential $\sum_{ab} Q_a^{u} u^b v_{ab}$,
with all numbers $Q_a$ being different, lies in the kernel of $\D$
defined in~\eqref{diff oper}. Then for every $a$, the quasi-exponential
$Q_a^{u} \sum_b u^b v_{ab}$ lies in the kernel of $\D$.
\end{lemma}

The lemma follows from the fact that exponential functions with dif\/ferent
exponents are linearly independent over the f\/ield of rational functions in $u$.

\section{Generating operator of the \XXX/ transfer matrices}
\label{Sec Generating}

\subsection[Yangian $Y(\glN)$]{Yangian $\boldsymbol{Y(\glN)}$}


{\bf 3.1.1.}
Let $e_{ab}$, $ a,b=1,\dots,N$, be the standard generators of
the complex Lie algebra $\glN$.
We have $\glN=\n^+\oplus\h\oplus\n^-$, where
\[
\n^+\>= \oplus _{a<b}\ \C\cdot e_{ab},
\qquad
\h = \oplus_{a=1}^N \ \C \cdot e_{aa},
\qquad
\n^-\>= \oplus _{a>b}\ \C \cdot e_{ab}.
\]
For an integral dominant $\glN$-weight $\La\in \h^*$,
denote by $M_\La$ the irreducible
f\/inite dimensional $\glN$-module with highest weight $\La$.

For a $\glN$-module $M$ and a weight $\mu\in\h^*$, denote by $M[\mu]\subset M$
the vector subspace of vectors of weight $\mu$.


\medskip

\noindent
{\bf 3.1.2.}
The Yangian $Y(\glN)$ is the unital associative algebra with
generators $T^{\{s\}}_{ab}$, $a,b=1,\dots,N$ and $s=1,2,\dots$.
Let
\[
T_{ab} = \delta_{ab} + \sum_{s=1}^\infty\, T^{\{s\}}_{ab} u^{-s} ,
\qquad
a,b=1,\dots,N .
\]
The def\/ining relations in $\Y$ have the form
\[
(u-v)\, [T_{ab}(u),T_{cd}(v)] = T_{cb}(v)T_{ad}(u) - T_{cb}(u)T_{ad}(v) ,
\]
for all $a$, $b$, $c$, $d$. The Yangian is a Hopf algebra with coproduct
\[
\Delta\ :\ T_{ab}(u) \ \mapsto \
\sum_{c=1}^N T_{cb}(u)\otimes T_{ac}(u)
\]
for all $a$, $b$.


\medskip

\noindent
{\bf 3.1.3.}
We identify the elements of $\End \,(\C^N)$ with $N\times N$-matrices.
Let $E_{ab}\in \End \,(\C^N)$ denote the matrix with the only nonzero entry $1$
at the intersection of the $a$-th row and $b$-th column.

Let $P = \sum_{a,b} E_{ab}\otimes E_{ba}$,
$R(u) = u + P \in \End\, (\C^N\otimes \C^N)$ and
$T(u) = \sum_{a,b} E_{ab}\otimes T_{ab}(u) \in
\End\,(\C^N) \otimes Y(\glN)((u^{-1}))$. Then the def\/ining relations for
the Yangian can be written as the following equation of series in $u^{-1}$
with coef\/f\/icients in $\End\,(\C^N) \otimes \End\,(\C^N)\otimes Y(\glN)$,
\[
R^{(12)}(u-v)\, T^{(13)}(u)\, T^{(23)}(v) =
T^{(23)}(v) \,T^{(13)}(u)\, R^{(12)}(u-v) .
\]


\medskip

\noindent
{\bf 3.1.4.}
A series $f(u)$ in $u^{-1}$ is called monic if $f(u)=1+O(u^{-1})$.

For a monic series $f(u)$, there is
an automorphism
\[
\chi_f\ :\ \Y\ \to\ \Y ,
\qquad
T(u)\ \mapsto\ f(u)\,T(u) .
\]
The f\/ixed point subalgebra in $\Y$ with respect to all automorphisms $\chi_f$
is called the Yangian $\Ys$. Denote by $Z\Y$ the center of the Yangian $\Y$.

\begin{proposition}
\label{propcenter}
The Yangian $\Y$ is isomorphic to the tensor product ${\Ys\otimes Z\Y}$.
\end{proposition}

See \cite{MNO} for a proof.

\medskip

\noindent
{\bf 3.1.5.}
Let $V$ be an irreducible f\/inite-dimensional $\Y$-module. Then there
exists a unique vector $v \in V$ such that
\begin{alignat*}{3}
& T_{ab}(u)\,v = 0, && a > b ,&
\\
& T_{aa}(u)\,v = c_a(u)\,v , \qquad && a=1,\dots,N ,&
\end{alignat*}
for suitable monic series $c_a(u)$. Moreover,
\begin{gather}
\frac{c_a(u)}{c_{a+1}(u)} = \frac{P_a(u+a)}{P_a(u+a-1)} ,
\qquad a=1,\dots,N-1,
\end{gather}
for certain monic polynomials $P_a(u)$.

The polynomials $P_1,\dots,P_{N-1}$ are called {\it the Drinfeld polynomials\/}
of the module $V$. The vector $v$ is called {\it a highest weight vector\/} and
the series $c_1(u),\dots,c_N(u)$ -- {\it the Yangian highest weights\/} of
the module $V$.

For any collection of monic polynomials $P_1,\dots, P_{N-1}$ there exists
an irreducible f\/inite-dimensional $\Y$-module $V$ such that the polynomials
$P_1,\dots, P_{N-1}$ are the Drinfeld polynomials of $V$. The module $V$ is
uniquely determined up to twisting by an automorphism of the form $\chi_f$.

The claim follows from Drinfeld's description of irreducible f\/inite-dimensional
$\Ys$-modu\-les~\cite{D} and Proposition~\ref{propcenter}.


\medskip

\noindent
{\bf 3.1.6.}
Let $V_1$, $V_2$ be irreducible f\/inite-dimensional $\Y$-modules with respective
highest weight vectors $v_1$, $v_2$. Then for the $\Y$-module $V_1\otimes V_2$,
we have
\begin{alignat*}{3}
& T_{ab}(u)\, v_1\otimes v_2 = 0 , && a > b ,&
\\
& T_{aa}(u)\, v_1\otimes v_2 = c_a^{(1)}(u)\,c_a^{(2)}(u)\,v_1\otimes v_2 ,
\qquad && a=1,\dots,N .&
\end{alignat*}

Let $W$ be the irreducible subquotient of $V_1\otimes V_2$ generated by
the vector $v_1\otimes v_2$. Then the Drinfeld polynomials of the module $W$
equal the products of the respective Drinfeld polynomials of the modules
$V_1$ and $V_2$.


\medskip

\noindent
{\bf 3.1.7.}
A f\/inite-dimensional irreducible
$\Y$-module $V$ will be called {\it polynomial} if{\samepage
\begin{gather}
\label{cN}
c_N(u) = \frac{P_N(u+N)}{P_N(u+N-1)} .
\end{gather}
for some monic polynomial $P_N(u)$.}

For any collection of monic polynomials $P_1,\dots, P_N$ there exists
a unique polynomial irreducible f\/inite-dimensional $\Y$-module $V$ such that
the polynomials $P_1,\dots, P_{N-1}$ are the Drinfeld polynomials of $V$ and
\Ref{cN} holds.


\medskip

\noindent
{\bf 3.1.8.}
There is a one-parameter family of automorphisms
\[
\rho_z \ :\ \Y \to \Y,\qquad T_{ab}(u) \mapsto T_{ab}(u-z) ,
\]
where in the right hand side,
$(u-z)^{-1}$ has to be expanded as a power series in $u^{-1}$.

The Yangian $\Y$ contains the universal enveloping algebra $U(\glN)$
as a Hopf subalgebra. The embedding is given by the formula
$e_{ab}\mapsto T^{\{1\}}_{ba}$ for all $a$, $b$.
We identify $U(\glN)$ with its image.

The evaluation homomorphism $\epsilon : \Y \to U(\glN)$ is def\/ined by the rule:
$ T^{\{1\}}_{ab} \mapsto e_{ba}$ for all $a$, $b$ and $T^{\{s\}}_{ab} \mapsto 0$
for all $a$, $b$ and all $s>1$.

For a $\glN$-module $V$ denote by $V(z)$ the $\Y$-module induced from
$V$ by the homomorphism $\epsilon \cdot \rho_z$. The module $V(z)$
is called the evaluation module with the evaluation point $z$.

Let $\bs \La = (\La^{(1)},\dots,\La^{(n)})$ be a collection of
integral dominant $\glN$-weights,
where
$\La^{(i)}=(\La_1^{(i)},\dots, \La_N^{(i)})$ for $i=1,\dots, n$.
For generic complex numbers $z_1,\dots,z_n$,
the tensor product of evaluation modules
\[
M_{\bs \La}(\bs z) = M_{\La^{(1)}}(z_1)\otimes\dots\otimes M_{\La^{(n)}}(z_n)
\]
is an irreducible $\Y$-module and the corresponding highest weight series
$c_1(u),\dots,$ $c_N(u)$ have the form
\begin{gather}
\label{cau}
c_a(u) = \prod_{i=1}^n\,\frac{u-z_i + \La_a^{(i)}}{u-z_i } .
\end{gather}
The corresponding Drinfeld polynomials are
\[
P_a(u)=\prod_{i=1}^n
\prod_{s=\La_{a+1}^{(i)}+1}^{\La_a^{(i)}}(u-z_i+s-a)
\]
for $a=1,\dots,N-1$.
The $\Y$-module $\MLz$ is polynomial if $\La_N^{(i)}\in\Z_{\geq 0}$
for all $i$.
Then the polynomial $P_N(u)$ has the form
\[
P_N(u)=\prod_{i=1}^n \prod_{s=1}^{\La_N^{(i)}}\,(u-z_i+s-N).
\]



\noindent
{\bf 3.1.9.}
Consider $\C^N$ as the $\glN$-module with highest weight $(1,0,\dots,0)$.

For any complex numbers $z_1,\dots,z_n$, all irreducible subquotients of
the $\Y$-module\\ ${\C^N(z_1)\otimes\dots\otimes \C^N(z_n)}$ are polynomial
$\Y$-modules. Moreover, for any polynomial irreducible f\/inite-dimensional
$\Y$-module $V$, there exist complex numbers $z_1,\dots,z_n$ such that $V$
is isomorphic to a subquotient of the $\Y$-module
$\C^N(z_1)\otimes\dots\otimes \C^N(z_n)$.

The numbers $z_1,\dots,z_n$ are determined by the formula
\[
\prod_{i=1}^n\,(u-z_i) = \prod_{a=1}^N \prod_{s=0}^{a-1} P_a(u+s) .
\]
This formula follows from consideration of the action of the center of
the Yangian $\Y$ in the module $\C^N(z_1)\otimes \dots \otimes \C^N(z_n)$.


\medskip

\noindent
{\bf 3.1.10.}
A f\/inite-dimensional $\Y$-module will be called {\it polynomial} if it is
the direct sum of tensor products of polynomial irreducible f\/inite-dimensional
$\Y$-modules.

If $V$ is a polynomial f\/inite-dimensional $\Y$-module, then for any
$a,b=1,\ldots,N$, the series $T_{ab}(u)|_V$ converges to an $\End\,(V)$-valued
rational function in $u$.

\medskip

\noindent
{\bf 3.1.11.}
Let ${\pi:U(\glN)\to\End\,(\C^N)}$ be the representation homomorphism for
the $\glN$-modu\-le~$\C^N$. Clearly, for any $x\in U(\glN)$ we have
\begin{gather}
\label{adjoint}
\bigl[\pi(x)\otimes 1+1\otimes x,\,T(u)\bigr] = 0.
\end{gather}

For a non-degenerate matrix $A\in\End\,(\C^N)$,
def\/ine an automorphism $\nu_A$ of $\Y$ by the formula
\[
(\id\otimes \nu_A)\bigl(T(u)\bigr)=
\sum_{ab}\ A^{-1} E_{ab}\>A\otimes T_{ab}(u).
\]

Let $V$ be a f\/inite-dimensional Yangian module with the representation
$\mu:Y(\glN)\to\End\,(V)$ and $\tilde\mu:GL_N\to\End\,(V)$ the corresponding
representation of the group $GL_N$. The automorphism $\nu_A$ induces a new
Yangian module structure $V^A$ on the same vector space with the representation
$\mu_A=\mu\circ\nu_A$. Formula \Ref{adjoint} yields that for any $x\in\Y$,
\begin{gather}
\label{muA}
\mu_A(x)=\tilde\mu(A)\,\mu(x)\bigl(\tilde\mu(A)\bigr)^{-1},
\end{gather}
that is, the $\Y$-modules $V$ and $V^A$ are isomorphic. In particular,
if $V$ is a polynomial irreducible f\/inite-dimensional $\Y$-module, then
$V^A$ is a polynomial irreducible f\/inite-di\-men\-sio\-nal $\Y$-module too.

\subsection[Universal difference operator]{Universal dif\/ference operator}


{\bf 3.2.1.}
Let $Q=(Q_{ab})$ be an $N\times N$-matrix. Def\/ine
\begin{gather}
\label{XabQ}
\Xe_{ab}(u, \tau)=
\delta_{ab}-\sum_{c=1}^N\,Q_{ac}\,T_{cb}(u)\,\tau^{-1},
\qquad a,b=1,\dots,N.
\end{gather}
If $V$ is a polynomial f\/inite-dimensional
$\Y$-module, then $\Xe_{ab}(u,\tau)$ acts on $V$-valued functions
in $u$,
\[
f(u)\ \mapsto\ \delta_{ab}f(u)-\sum_{c=1}^N\,Q_{ac}\,T_{cb}(u)f(u-1) .
\]
Following \cite{T}, introduce the dif\/ference operator
\begin{gather}
\label{Dg}
\D(u,\tau)=\sum_{\si\in S_N} (-1)^\si\,
\Xe_{1\,\si_1}(u,\tau)\,\Xe_{2\,\si_2}(u,\tau)\cdots
\Xe_{N\si_N}(u,\tau),
\end{gather}
where the sum is over all permutations $\si$ of $\{1,\dots,N\}$.
The operator $\D(u,\tau)$ will be called {\it the universal difference
operator\/} associated with the matrix $Q$.{\samepage

\begin{lemma}
\label{rowdet}
Let $\pi$ be a map $\{1,\dots,N\}\to\{1,\dots,N\}$.
If $\pi$ is a permutation of $\{1,\dots,N\}$, then
\[
\sum_{\si\in S_N} (-1)^\si\,
\Xe_{\pi_1\si_1}(u,\tau)\,\Xe_{\pi_2\si_2}(u,\tau)\cdots
\Xe_{\pi_N\si_N}(u,\tau) = (-1)^\pi\,\D(u,\tau) .
\]
If $\pi$ is not bijective, then
\[
\sum_{\si\in S_N} (-1)^\si\,
\Xe_{\pi_1\si_1}(u,\tau)\,\Xe_{\pi_2\si_2}(u,\tau)\cdots
\Xe_{\pi_N\si_N}(u,\tau) = 0 .
\]
\end{lemma}

The statement is Proposition~4.10 in~\cite{MTV2}.

}

\medskip

\noindent
{\bf 3.2.2.}
Introduce the coef\/f\/icients $\Te_0(u),\dots,\Te_N(u)$ of $\D(u,\tau)$:
\[
\D(u,\tau)=\sum_{k=0}^N\,(-1)^k\,\Te_k(u)\,\tau^{-k}.
\]
The coef\/f\/icients $\Te_k(u)$ are called {\it the transfer matrices of
the \XXX/ type model} associated with~$Q$.
The main properties of the transfer matrices:
\begin{enumerate}\itemsep=0pt
\item[(i)]
the transfer matrices commute: $[\Te_k(u),\Te_l(v)] = 0$
for all $k$, $l$, $u$, $v$,
\item[(ii)]
if $Q$ is a diagonal matrix, then the transfer matrices preserve
the $\glN$-weight: $\![\Te_k(u),\! e_{aa}]{=}0$ for all $k$, $a$, $u$,
\item[(iii)]
if $Q$ is the identity matrix, then the transfer matrices commute with the
subalgebra $U(\glN)$: $[\Te_k(u), x] = 0$ for all $k,u$ and $x\in U(\glN)$,
\end{enumerate}
see \cite{T, MTV2}.


\medskip

\noindent
{\bf 3.2.3.}
Evidently, $\Te_0(u)=1$.
We also have\ $\Te_N(u) = \det Q \, {\rm qdet}\, T(u)$, where
\[
{\rm qdet}\, T(u) =
\sum_{\sigma\in S_N}
(-1)^\si\,
T_{1\,\si_1}(u)\,T_{2\,\si_2}(u-1)\cdots
T_{N\si_N}(u-N+1),
\]
and ${\rm qdet}\, T(u) = 1 + O(u^{-1})$.

\begin{theorem}
The coefficients of the series ${\rm qdet}\,T(u)$ are free generators of
the center of the Yangian $\Y$.
\end{theorem}

See \cite{MNO} for a proof.


\medskip

\noindent
{\bf 3.2.4.}
If $V$ is a polynomial f\/inite-dimensional
$\Y$-module, then the universal operator $\D(u,\tau)$
induces a dif\/ference operator acting on $V$-valued functions in
$u$. This operator
will be called {\it the universal difference
operator} associated with $Q$ and $V$ and denoted by
$\D_{Q,V}(u,\tau)$. The linear
operators $\Te_k(u)|_{V} \in \End\,(V)$ will be
called {\it the transfer matrices} associated with $Q$ and
$V$ and denoted by $\Te_{k,Q,V}(u)$.
They are rational functions in $u$.

\begin{example}
Let
\[
V=M_{\bs \La}(\bs z)=M_{\La^{(1)}}(z_1)\otimes\dots\otimes
M_{\La^{(n)}}(z_n),
\]
Consider the algebra $\bigl(U(\glN)\bigr)^{\otimes n}$.
For $a,b=1,\dots,N$ and $i=1,\dots,n$, def\/ine
\[
e_{ab}^{(i)} = 1^{\otimes(i-1)}\otimes e_{ab}\otimes 1^{\otimes(n-i)},
\qquad
L_{ab}^{(i)}(u,z) = \delta_{ab}+\frac 1{u-z}\,e_{ba}^{(i)}.
\]
The operator $T_{ab}(u)$ acts on $\MLz$ as
\[
\sum_{c_1,\dots, c_{n-1}=1}^{N} \!\! L_{a\,c_{n-1}}^{(n)}(u-z_n)\,
L_{c_{n-1}\,c_{n-2}}^{(n-1)}(u-z_{n-1})\cdots L_{c_2\,c_1}^{(2)}(u-z_2)\,
L_{c_1b}^{(1)}(u-z_1).
\]
\end{example}

\begin{lemma}
\label{sec on nu-A 2}
If $\D_Q(u,\tau)$ is the universal difference operator associated with
the matrix $Q$ and $\nu_A : Y(\glN) \to Y(\glN)$ is the automorphism defined
in Section {\rm 3.1.11}, then
\begin{gather}
\label{nuAD}
\nu_A(\D_Q(u,\tau)) = \D_{AQA^{-1}}(u,\tau)
\end{gather}
is the universal difference operator associated with the matrix $AQA^{-1}$.
\end{lemma}

\begin{proof}
Consider matrices $Q=(Q_{ab})$, $T=(T_{ab}(u))$, and
$\Xe=(\Xe_{ab}(u,\tau))$. Then
formula \Ref{XabQ}, def\/ining $\Xe_{ab}$, may be
read as $\Xe = 1 - Q T \tau^{-1}$.
Formula \Ref{Dg} for the universal dif\/ference operator may be understood
as the row determinant of the matrix $\Xe$, $\D(u,\tau) = \det\Xe$.
The def\/inition of the automorphism $\nu_A$ reads as
$(\id\otimes\nu_A)(T) = A^{-1} T A.$
Then we have
\[
\nu_A(\D(u,\tau)) = \det\big( 1 - Q A^{-1} T A \tau^{-1} \big) =
\det\big( A^{-1} \big( 1 - A Q A^{-1} T \tau^{-1} \big) A \big).
\]
Now formula \Ref{nuAD} will be proved
if we were able to write that last determinant as the product:
\[
\det\big( A^{-1} \big( 1 - A Q A^{-1} T \tau^{-1} \big) A \big) =
\det\big(A^{-1}\big) \det \big( 1 - A Q A^{-1} T \tau^{-1} \big) \det A .
\]
The last formula may be proved the same way as the standard formula
$\det MN = \det M\,\det N$ in ordinary linear algebra, using two
observations. The f\/irst is that the entries of $A$ are numbers and commute with
the entries of $\Xe$. The second observation is Lemma~\ref{rowdet} describing
the transformations of the row determinant of $\Xe$ with respect to row
replacements.
\end{proof}



\noindent
{\bf 3.2.5.}
Let $V$ be a polynomial f\/inite-dimensional Yangian module and
$\tilde\mu:G_N\to GL(V)$ the associated $GL_N$-representation.
Then formulae \Ref{muA} and \Ref{nuAD} yield
\[
\D_{AQA^{-1}}(u,\tau)|_V =
\tilde\mu(A)\,\D_Q(u,\tau)|_{V}\,\tilde \mu\big(A^{-1}\big) .
\]

\subsection{More properties of transfer matrices}
\label{Properties of transfer matrices}

Let $Q = \diag\,(Q_1,\dots,Q_N)\, \in GL_N$.
Let $\bs \La = (\La^{(1)},\dots,\La^{(n)})$ be a collection of
integral domi\-nant~$\glN$-weights,
where $\La^{(i)}=(\La_1^{(i)},\dots, \La_N^{(i)})$ and
$\La_N^{(i)}=0$ for $i=1,\dots, n$.

For $\bs m = (m_1,\dots,m_N)$ denote by
$M_{\bs \La}(\bs z)[\bs m] \,\subset\, M_{\bs \La}(\bs z)$ the
weight-subspace of $\glN$-weight~$\bs m$ and by
$\sing\, M_{\bs \La}(\bs z)[\bs m]
\subset M_{\bs \La}(\bs z)[\bs m]$ the subspace of $\glN$-singular vectors.

Note that the subspace $\sing\, M_{\bs \La}(\bs z)[\bs m]$ may be nonzero only
if $\bs m $ is dominant integral, i.e.\ the integers $m_1,\dots, m_N$ have to
satisfy the inequalities $m_1\geq \dots \geq m_N$.

Consider the universal dif\/ference operator associated with $Q$ and
$M_{\bs \La}(\bs z)$. Then the associated transfer matrices preserve
$M_{\bs \La}(\bs z)[\bs m]$ and we may consider the universal dif\/ference
operator $\D_{Q,\,M_{\bs \La}(\bs z)[\bs m]}(u,\tau)$ acting on
$M_{\bs \La}(\bs z)[\bs m]$-valued functions of $u$. We may write
\[
\D_{Q,\,M_{\bs \La}(\bs z)[\bs m]}(u,\tau) = \sum_{k=0}^N\,(-1)^{k}\,
\Te_{k,Q,\,M_{\bs \La}(\bs z)[\bs m]}(u)\,\tau^{-k} .
\]
As we know $\Te_{0,Q,M_{\bs \La}(\bs z)[\bs m]}(u)=1$, and we have
\[
\Te_{k,Q,M_{\bs \La}(\bs z)[\bs m]}(u) =
\Te_{k0} + \Te_{k1}u^{-1} + \Te_{k2}u^{-2} + \cdots
\]
for suitable $\Te_{ki} \in \End\,(M_{\bs \La}(\bs z)[\bs m])$.

\begin{theorem}
\label{thm on coefficients}
The followings statements hold.
\begin{enumerate}\itemsep=0pt
\item[(i)]
The operators $\Te_{10},\Te_{20},\dots,\Te_{N0}$ and
$\Te_{11},\Te_{21},\dots,\Te_{N1}$ are scalar operators.
Moreover, the following relations hold:
\begin{gather*}
x^N+\,\sum_{k=1}^N \,(-1)^{k}\,\Te_{k0}\,x^{N-k} = \prod_{i=1}^N\,(x-Q_i) ,
\\
\sum_{k=1}^N \,(-1)^{k}\,\Te_{k1}\,x^{N-k} =
- \prod_{i=1}^N\,(x-Q_i) \sum_{j=1}^N\,\frac{m_j\,Q_j}{x-Q_j} .
\end{gather*}

\item[(ii)]
For $k=1,\dots,N-1$,
\begin{gather*}
\Te_{N,Q,M_{\bs \La}(\bs z)[\bs m]}(u) =
\left(\prod_{i=1}^N Q_i\right) \prod_{s=1}^n
\prod_{i=1}^{N-1}\,\frac{u-z_s+\La^{(s)}_i-i+1}{u-z_s-i+1} ,
\\
\Te_{k,Q,M_{\bs \La}(\bs z)[\bs m]}(u) =
\tilde{\Te}_{k}(u)
\prod_{i=1}^k\,\prod_{s=1}^n\,\frac1{u-z_s-i+1} ,
\end{gather*}
where $\tilde{\Te}_{k}(u)$ is a polynomial in $u$ of degree $nk$.
\end{enumerate}
\end{theorem}

\begin{proof}
Part (i) follows from Proposition B.1 in~\cite{MTV2}. Part (ii) follows
from the def\/inition of the universal dif\/ference operator and the fact that
the coef\/f\/icients of the series $\Te_N(u)$ belong to the center of $\Y$.
\end{proof}



\noindent
{\bf 3.3.1.}
Assume that $Q$ is the identity matrix. Then the associated transfer matrices
preserve $\sing\,M_{\bs \La}(\bs z)[\bs m]$ and we may consider the universal
dif\/ference operator $\D_{Q=1,\,\sing M_{\bs \La}(\bs z)[\bs m]}(u,\tau)$ acting
on $\sing\,M_{\bs \La}(\bs z)[\bs m]$-valued functions of $u$. We may write
\begin{gather}
\label{formula for Se}
\D_{Q=1,\,\sing M_{\bs \La}(\bs z)[\bs m]}(u,\tau)\,\tau^N
= \sum_{k=0}^N \, (-1)^k\, {\Se}_{k}(u) \,(\tau - 1)^{N-k}
\end{gather}
for suitable coef\/f\/icients ${\Se}_{k}(u)$.

Note that the operators $\Se_k(u)$ coincide with the action in
$\sing\,M_{\bs \La}(\bs z)[\bs m]$ of the modif\/ied transfer matrices
$\Tee_k(u)$ from formula~(10.4) of \cite{MTV2}.

\begin{theorem}
\label{thm Q=1}
The following three statements hold.
\begin{enumerate}\itemsep=0pt
\item[(i)]
We have ${\Se}_0(u) = 1$.
\item[(ii)]
For $k=1,\dots, N$, the coefficient $ {\Se}_{k}(u)$ has the following
Laurent expansion at $u=\infty$:
\[
{\Se}_{k}(u) = {\Se}_{k,0} u^{-k} + {\Se}_{k,1} u^{-k-1} + \cdots ,
\]
where the operators ${\Se}_{1,0},\dots,{\Se}_{N,0}$ are scalar operators.
\item[(iii)]
For all $d$ we have
\[
\sum_{k=0}^N (-1)^k\,\Se_{k,0}\prod_{j=0}^{N-k-1} (d-j)
= \prod_{s=1}^N\,(d - m_s - N + s) .
\]
\end{enumerate}
\end{theorem}

\begin{proof}
Part (i) is evident.

Since $Q$ is the identity matrix, formula \Ref{XabQ} reads now as follows
\begin{gather}
\label{Xab}
\Xe_{ab}(u) =\delta_{ab}-T_{ab}(u)\,\tau^{-1}=
\bigl(\delta_{ab}\,(\tau-1)-T_{ab}^{\{1\}}\,u^{-1}-
O\big(u^{-2}\big)\bigr)\,\tau^{-1}.
\end{gather}
Then part (ii) is straightforward from formulas \Ref{formula for Se}
and~\Ref{Dg}.

Let $v$ be any vector in $\sing\,M_{\bs \La}(\bs z)[\bs m]$ and $d$ any number.
To prove part (iii) we apply the dif\/ference operators in formula
\Ref{formula for Se} to the function $v u^d$. The expansion at inf\/inity
of result of the application of the right side is
\[
u^{d-N}\left(\sum_{k=0}^N (-1)^k \left[\prod_{j=0}^{N-k-1} (d-j)\right]
\Se_{k,0}\,v + O\big(u^{-1}\big)\right).
\]
So it remains to show that
\[
\D_{Q=1,\sing M_{\bs \La}(\bs z)[\bs m]}(u,\tau)\,\tau^N v\,u^d=u^{d-N}
\left(\left[\prod_{s=1}^N (d - m_s - N +s)\right] v + O\big(u^{-1}\big)\right),
\]
as $u$ goes to inf\/inity.
Since $\tau^N u^d=u^d\bigl(1+O(u^{-1})\bigr)$, the last formula is equivalent to
\begin{gather}
\label{formula leading term}
\D_{Q=1,\,\sing M_{\bs \La}(\bs z)[\bs m]}(u,\tau)\,v\,u^d =u^{d-N}
\left(\left[\prod_{s=1}^N(d - m_s - N + s)\right] v + O\big(u^{-1}\big)\right).
\end{gather}
To prove \Ref{formula leading term}, observe that according to formula
\Ref{Xab} we have
\[
\Xe_{ab}(u)\, v u^d=u^{d-1}\,\Bigl(
\bigl(d\,\delta_{ab}-T_{ab}^{\{1\}}(u)\bigr)\,v+ O\big(u^{-1}\big)\Bigr).
\]
Applying this remark to formula \Ref{Dg} we get
\begin{gather*}
\D_{Q=1,\,\sing M_{\bs \La}(\bs z)[\bs m]}(u,\tau)\,v\,u^d
\\
\qquad {}=u^{d-N}\left(\,\sum_{\si\in S_N} (-1)^\si\,
\big((d-N+1)\,\delta_{1,\si_1} - T^{\{1\}}_{1\si_1}\big)\cdots
\big(d \,\delta_{N,\si_N} - T^{\{1\}}_{N\si_N}\big)\,v+ O\big(u^{-1}\big)\right).
\end{gather*}
Each element $T^{\{1\}}_{ab}$ acts as $\sum\limits_{i=1}^n e^{(i)}_{ba}$
in $M_{\bs\La}(\bs z)$, which corresponds to the standard $\glN$-action
in the $\glN$-module $M_{\bs\La}$. Since $v$ is a $\glN$-singular vector
of a $\glN$-weight $(m_1,\dots, m_N)$, we have
\[
\sum_{\si\in S_N} (-1)^\si\big((d-N+1)\,\delta_{1,\si_1} - T^{\{1\}}_{1\si_1}\big)\cdots
\big(d \,\delta_{N,\si_N} - T^{\{1\}}_{N\si_N}\big)v=
\left[\prod_{s=1}^N (d - m_s - N + s)\right] v .
\]
Indeed only the identity permutation contributes nontrivially to the sum
in the left side. This proves part (iii) of Theorem \ref{thm Q=1}.
\end{proof}

\subsection{First main result}
\label{The first main}
\begin{theorem}
\label{The first main thm}
Let $V$ be a polynomial finite-dimensional $\Y$-module and $Q\in GL_N$.
Consider the universal difference operator $\D_{Q,V}(u,\tau)$ associated
with $Q$ and $V$. Then the kernel of $\D_{Q,V}(u,\tau)$ is generated by
quasi-exponentials.
\end{theorem}

The theorem will be proved in Section \ref{Bethe ansatz}.

\section[Continuity principle for difference operators\\ with
quasi-exponential kernel]{Continuity principle for dif\/ference operators\\ with
quasi-exponential kernel}
\label{Sec Cont}

\subsection{Independent quasi-exponentials}
\label{sec Quasi-exponentials}
Let $V$ be a complex vector space of dimension $d$. Let $p \in \C[u]$
be a monic polynomial of degree~$k$.
Consider the dif\/ferential equation
\[
p\<\left(\frac {d}{du}\right)\< f(u) = 0
\]
for a $V$-valued function $f(u)$.
Denote by $W_p$ the complex vector space of its solutions. The map{\samepage
\[
\delta_p\ :\ W_p\ \to\ V^{\oplus k} , \quad w \ \mapsto\ \big(w(0),w'(0),\dots,
w^{(k-1)}(0)\big),
\]
assigning to a solution its initial condition at $u=0$, is an isomorphism.}

Let $\lambda_1,\dots,\lambda_l$ be all distinct roots of the polynomial $p$ of
multiplicities $k_1,\dots,k_l$, respectively. Let $v_1,\dots,v_d$ be a basis
of $V$. Then the quasi-exponentials
\begin{gather}
\label{basis}
e^{\lambda_j u} u^a v_b ,\qquad j=1,\dots,l,\quad
a=0,\dots,k_j-1 ,\quad b = 1,\dots, d ,
\end{gather}
form a basis in $W_p$.

\begin{lemma}
\label{lem on lin indep}
Assume that $\lambda_a-\lambda_b \notin 2\pi i\,\Z$ for all $a\neq b$.
Then the $kd$ functions listed in \Ref{basis} are linear independent over
the field of one-periodic functions.
\end{lemma}

\subsection[Admissible difference operators]{Admissible dif\/ference operators}
Let
$A_0(u),\dots, A_N(u)$ be $\End\,(V)$-valued rational functions in $u$.
Assume that each of these functions has limit as $u \to \infty$ and
$A_0(u)=1$ in $\End\,(V)$. For every $k$, let
\[
A_k(u) = A^\infty_{k,0} + A^\infty_{k,1}
u^{-1} + A^\infty_{k,2} u^{-2} + \cdots
\]
be the Laurent expansion at inf\/inity.
Consider the algebraic equation
\begin{gather}
\label{char eqn}
\det \left( A_{N,0} + x A_{N-1,0} + \dots + x^{N-1} A_{1,0}
+ x^N \right) = 0
\end{gather}
with respect to variable $x$ and the dif\/ference operator
\[
\D = \sum_{k=0}^N\, A_k(u)\, \tau^{-k}
\]
acting on $V$-valued functions in $u$. Equation \Ref{char eqn} will be
called the characteristic equation for the dif\/ference operator $\D$.

The operator $\D$ will be called
{\it admissible at infinity\/} if $\det\, A_{N,0} \neq 0$, or equivalently,
if $x=0$ is not a root of the characteristic equation.

\begin{example}
Let $V$ be the tensor product of polynomial f\/inite-dimensional
irreducible $Y(\glN)$-modules and $Q\in GL_N$.
Let $\D_{Q,V}(u,\tau)$ be the associated
universal dif\/ference operator. Then $\D_{Q,V}(u,\tau)$ is admissible at
inf\/inity, see Sections 3.2.3 
and \ref{Properties of transfer matrices}.
\end{example}

\begin{lemma}
\label{lem of Q}
Assume that $\D$ is admissible at infinity and a nonzero $V$-valued
quasi-exponential $Q^u(u^dv_d + u^{d-1}v_{d-1} + \dots + v_0)$ lies
in the kernel of $\D$. Then $Q$ is a root of the characteristic
equation \Ref{char eqn}.
\end{lemma}

\subsection{Continuity principle}
\label{sec continuity principle}

Let
$A_0(u,\epsilon),\dots, A_N(u,\ep)$ be $\End\,(V)$-valued rational functions
in $u$ analytically depending on $\ep \in [0,1)$. Assume that
\begin{enumerate}\itemsep=0pt
\item[$\bullet$]
for every $\ep \in [0,1)$ the dif\/ference operator
$\D_\ep = \sum\limits_{k=0}^N A_k(u,\ep) \tau^{-k}$
is admissible at inf\/inity,
\item[$\bullet$]
for every $\ep \in (0,1)$ the kernel of $\D_\ep$ is
generated by quasi-exponentials,
\item[$\bullet$]
there exists a natural number $m$ such that for every $\ep\in(0,1)$ all
quasi-exponentials generating the kernel of $\D_\ep$ are of degree less than
$m$.
\end{enumerate}

\begin{theorem}
\label{thm continuity principle}
Under these conditions the kernel of the difference operator $\D_{\ep=0}$ is
generated by quasi-exponentials.
\end{theorem}
\begin{proof}
For every $\ep \in [0,1)$, the characteristic equation for $\D_\ep$
has $Nd$ roots counted with multiplicities.
As $\ep$ tends to $0$ the roots of the characteristic equation
of $\D_\ep$ tend
to the roots of the characteristic equation of $\D_{\ep=0}$.
All these roots are nonzero numbers.
For small positive $\ep$ the set of multiplicities of
roots does not depend on $\ep$.

\medskip

The following lemma is evident.

\begin{lemma}
There exist
\begin{enumerate}\itemsep=0pt
\item[$\bullet$]
a number $\bar \epsilon$ with $0<\bar \epsilon \leq 1$,
\item[$\bullet$]
for any $\ep$, $0<\epsilon<\bar \epsilon$,
a way to order the roots of the characteristic equation of $\D_\ep$
(we denote the ordered roots by $Q_1^\ep,\dots, Q_{Nd}^\ep$),
\item[$\bullet$]
a way to assign the logarithm $q_j^\ep$ to every root $Q_j^\ep$
\end{enumerate}
such that for every $j$ the number $q_j^\ep$ continuously depends
on $\ep$ and $q_j^\ep = q_l^\ep$ whenever $Q_j^\ep = Q_l^\ep$.
\end{lemma}

Let $m$ be the number described in Section~\ref{sec continuity principle}.
For every $\ep$, $0<\epsilon<\bar \epsilon$, we
def\/ine $p_\ep \in \C[u]$ to
be the monic polynomial of degree $k=mNd$,
whose set of roots consists of $m$ copies of each of
the numbers $q_1^\ep,\dots, q_{Nd}^\ep$.

Let $W_{p_\ep}$ be the $kd$-dimensional vector space of
quasi-exponentials assigned to the polynomial~$p_\ep$
in Section \ref{sec Quasi-exponentials}.
By assumptions of Theorem \ref{sec Quasi-exponentials}, for every
$\ep$, $0<\epsilon<\bar \epsilon$,
the space $W_{p_\ep}$ contains an $Nd$-dimensional subspace $U_\ep$
generating the kernel of $\D_\ep$.
This subspace determines a point in the
Grassmannian $Gr(W_{p_\ep}, Nd)$ of $Nd$-dimensional subspaces of
$W_{p_\ep}$.

The map $\rho_{p_\ep} : W_{p_\ep} \to V^{\oplus k}$ identif\/ies
the Grassmannian ${\rm Gr}(W_{p_\ep}, Nd)$ with the Grassmannian
${\rm Gr}(V^{\oplus k}, Nd)$. The points $\rho_{p_\ep} (U_\ep)$ all lie in the
compact manifold ${\rm Gr}(V^{\oplus k}, Nd)$ and the set of all such points
has an accumulation point $\tilde U \in
{\rm Gr}(V^{\oplus k}, Nd)$ as $\ep$ tends to zero.
Then $\rho_{p_{\ep=0}}^{-1}(\tilde U)\subset
{\rm Gr}(W_{p_{\ep=0}}, Nd)$ is an $Nd$-dimensional
subspace of $V$-valued quasi-exponentials. Using Lemma \ref{lem on lin indep},
we conclude that the space
$\rho_{p_{\ep=0}}^{-1}(\tilde U)$ generates the kernel of $\D_{\ep=0}$.
\end{proof}

\section {Bethe ansatz}
\label{Bethe ansatz}

\subsection{Preliminaries}
Consider $\C^N$ as the $\glN$-module with highest weight
$(1,0,\dots,0)$. For complex numbers $z_1,\dots,\alb gz_n$,
denote $\bs z = (z_1,\dots,z_n)$, and
\[
M(\bs z) = \C^N(z_1)\otimes \dots \otimes \C^N(z_n) ,
\]
which is a polynomial $\Y$-module.
Let
\[
M(\bs z) = \oplus_{m_1,\dots,m_N} M(\bs z)[m_1,\dots,m_N]
\]
be its $\glN$-weight decomposition with respect to the Cartan subalgebra of
diagonal matrices. The weight subspace
$M(\bs z)[m_1,\dots,m_N]$ is nonzero if and only if
$m_1,\dots, m_N $ are non-negative integers.

Assume that $Q = \diag\,(Q_1,\dots,Q_N)$ is a diagonal
non-degenerate $N\times N$-matrix with distinct diagonal entries,
and consider the universal dif\/ference operator
\[
\D_{Q,\,M(\bs z)}(u,\tau) = \sum_{k=0}^N
(-1)^k\,
\Te_{k,Q, M(\bs z)}(u)\,\tau^{-k}.
\]
associated with $M(\bs z)$ and $Q$.
Acting on $M(\bs z)$-valued functions the operator
$\D_{Q, M(\bs z)}(u,\tau)$ preserves the weight decomposition.

In this section we shall study the kernel of this dif\/ference operator,
restricted to $M(\bs z)[m_1,\dots,\alb m_N]$-valued functions.

\subsection{Bethe ansatz equations associated with a weight subspace}
\label{Bethe ansatz equations associated with a weight subspace}
Consider a nonzero weight subspace $M(\bs z)[m_1,\dots,m_N]$. Introduce
$\bs l = (l_1,\dots,l_{N-1})$ with $l_j \ =\ m_{j+1} + \dots + m_N$.
We have $n\geq l_1\geq \dots \geq l_{N-1}\geq 0$.
Set $l_0=l_N=0$ and $l=l_1 + \dots + l_{N-1}$.
We shall consider functions of $l$ variables
\[
\bs t = \big(t^{(1)}_{1},\dots, t^{(1)}_{l_1},
t^{(2)}_{1},\dots, t^{(2)}_{l_2},\dots,
t^{(N-1)}_{1},\dots, t^{(N-1)}_{l_{N-1}} \big) .
\]
The following system of $l$ algebraic equations with respect to $l$ variables
$\bs t$ is called {\it the Bethe ansatz equations} associated with
$M(\bs z)[m_1,\dots,m_N]$ and $Q$,
\begin{gather}
\label{BAE}
Q_1 \prod_{s=1}^n \big(t^{(1)}_j - z_s + 1\big)
\prod_{\satop{j'=1}{j'\neq j}}^{l_1}
\big(t^{(1)}_j - t^{(1)}_{j'} - 1 \big)
\prod_{j'=1}^{l_2} \big(t^{(1)}_j - t^{(2)}_{j'}\big)
\\
\qquad {}=
Q_2 \prod_{s=1}^n \big(t^{(1)}_j - z_s\big)\prod_{\satop{j'=1}{j'\neq j}}^{l_1}
\big(t^{(1)}_j - t^{(1)}_{j'} + 1 \big)
\prod_{j'=1}^{l_2} \big(t^{(1)}_j - t^{(2)}_{j'} - 1\big) ,
\notag
\\
Q_a \prod_{j'=1}^{l_{a-1}} \big(t^{(a)}_j - t^{(a-1)}_{j'} + 1\big)
\prod_{\satop{j'=1}{j'\neq j}}^{l_a} \big(t^{(a)}_j - t^{(a)}_{j'} - 1 \big)
\prod_{j'=1}^{l_{a+1}} \big(t^{(a)}_j - t^{(a+1)}_{j'}\big)
\notag
\\
\qquad{} =
Q_{a+1} \prod_{j'=1}^{l_{a-1}} \big(t^{(a)}_j - t^{(a-1)}_{j'}\big)
\prod_{\satop{j'=1}{j'\neq j}}^{l_a} \big(t^{(a)}_j - t^{(a)}_{j'} + 1 \big)
\prod_{j'=1}^{l_{a+1}} \big(t^{(a)}_j - t^{(a+1)}_{j'}-1\big) ,
\notag
\\
Q_{N-1} \prod_{j'=1}^{l_{N-2}} \big(t^{(N-1)}_j - t^{(N-2)}_{j'} + 1\big)
\prod_{\satop{j'=1}{j'\neq j}}^{l_{N-1}} \big(t^{(N-1)}_j - t^{(N-1)}_{j'} - 1 \big)
\notag
\\
\qquad{}=
Q_{N} \prod_{j'=1}^{l_{N-2}} \big(t^{(N-1)}_j - t^{(N-2)}_{j'}\big)
\prod_{\satop{j'=1}{j'\neq j}}^{l_{N-1}} \big(t^{(N-1)}_j - t^{(N-1)}_{j'} + 1 \big) .
\notag
\end{gather}
Here the equations of the f\/irst group are labeled by $j=1,\dots,l_1$,
the equations of the second group are labeled by $a=2,\dots,N-2$,
$j=1,\dots,l_a$, the equations of
the third group are labeled by $j=1,\dots,l_{N-1}$.

A solution $\tilde {\bs t}$ of system \Ref{BAE}
will be called {\it off-diagonal\/} if
\ ${\tilde t^{(a)}_{j} \neq \tilde t^{(a)}_{j'}}$ \,for any
$a=1,\dots,\alb N-1$, $1\leq j \leq j' \leq l_a$, and
$\tilde t^{(a)}_{j} \neq \tilde t^{(a+1)}_{j'}$
for any
$a=1,\dots,N-2$, $j= 1,\dots,l_a, \, j'= 1,\dots,l_{a+1}$.

\subsection{Weight function and Bethe ansatz theorem}
Denote by $\omega (\bs t,\bs z)$ {\it the universal weight function}
associated with the weight subspace $M(\bs z)[m_1,\alb\dots,m_N]$.
The universal weight function is def\/ined in formula~(6.2) in \cite{MTV2},
see explicit formula~\Ref{explicit} below. At this moment, it is enough for us
to know that this function is an $M(\bs z)[m_1,\dots,m_N]$-valued polynomial
in $\bs t$, $\bs z$.

If $\tilde{\bs t}$ is an of\/f-diagonal solution of the Bethe ansatz equations,
then the vector $\omega(\tilde {\bs t},\bs z) \in
M(\bs z)[m_1,\dots,m_N]$ is called {\it the Bethe vector} associated with
$\tilde{\bs t}$.{\samepage

\begin{theorem}
\label{thm on Bethe ansatz}
Let $Q$ be a diagonal matrix and $\tilde{\bs t}$ an off-diagonal solution
of the Bethe ansatz equations \Ref{BAE}. Assume that the Bethe vector
$\omega(\tilde {\bs t},\bs z)$ is nonzero. Then the Bethe vector is
an eigenvector of all transfer-matrices $\Te_{k,Q,M(\bs z)}(u)$, $k=0,\dots,N$.
\end{theorem}

The statement follows from Theorem~6.1 in~\cite{MTV2}.
For ${k=1}$, the result is established in~\cite{KR}.

}

\medskip

The eigenvalues of the Bethe vector are given by the following construction.
Set
\begin{gather*}
\chi^1(u, \bs t,\bs z) =
Q_1 \prod_{s=1}^n \frac{u-z_s + 1}{u-z_s}
\prod_{j=1}^{l_1} \frac{u- t^{(1)}_j - 1}{u-t^{(1)}_j} ,
\\
\chi^a(u, \bs t,\bs z) = Q_a \prod_{j=1}^{l_{a-1}}
\frac{u- t^{(a-1)}_j + 1}{u-t^{(a-1)}_j}\
\prod_{j=1}^{l_a} \frac{u-t^{(a)}_j - 1}{u-t^{(a)}_j} ,
\end{gather*}
for $ a=2,\dots, N$. Def\/ine the functions
$\lambda_k(u,\bs t,\bs z)$ by the formula
\[
\big(1-\chi^1(u,\bs t,\bs z)\,\tau^{-1}\big)\cdots
\big(1-\chi^N(u,\bs t,\bs z)\,\tau^{-1}\big) =
\sum_{k=0}^N \,(-1)^k\,\lambda_k(u,\bs t,\bs z)\,\tau^{-k} .
\]
Then
\[
\Te_{k, Q, M(\bs z)}(u)\,\omega(\tilde{\bs t},\bs z) =
\lambda_k(u,\tilde{\bs t},\bs z)\,
\omega(\tilde{\bs t},\bs z)
\]
for $k=0,\dots, N$, see Theorem 6.1 in \cite{MTV2}.

\subsection[Difference operator associated with an off-diagonal solution]{Dif\/ference operator associated with an of\/f-diagonal solution}

Let $\tilde{\bs t}$ be an of\/f-diagonal solution of the Bethe ansatz equations.
The scalar dif\/ference operator
\[
\D_{\tilde {\bs t}}(u,\tau) = \sum_{k=0}^N \,
(-1)^k\,\lambda_k(u,\tilde{\bs t},\bs z)\,\tau^{-k}
\]
will be called the associated {\it fundamental difference operator}.

\begin{theorem}
\label{thm on fund operator}
The kernel of $\Dt$ is generated by quasi-exponentials of degree
bounded from above by a function in $n$ and $N$.
\end{theorem}

This is Proposition~7.6 in \cite{MV3}, which is a generalization of
Proposition~4.8 in \cite{MV2}.

\subsection{Completeness of the Bethe ansatz}
\begin{theorem}
\label{thm on completeness}
Let $z_1,\dots,z_n$ and $Q=\diag\,(Q_1,\dots,Q_N)$ be generic.
Then the Bethe vectors form a basis in $M(\bs z)[m_1,\dots,m_N]$.
\end{theorem}

Theorem \ref{thm on completeness} will be proved in
Section~\ref{prf of thm on completeness}.

\begin{corollary}
Theorems {\rm \ref{thm on completeness}} and~{\rm \ref{thm continuity principle}}
imply Theorem~{\rm \ref{The first main thm}}.
\end{corollary}

\begin{proof}
Theorems \ref{thm on fund operator} and \ref{thm on completeness}
imply that the statement of Theorem \ref{The first main thm} holds
if the tensor product $M (\bs z)$ is considered for generic $\bs z$ and
generic diagonal $Q$. Then according to the remark in
Section 3.2.5, 
the statement of Theorem \ref{The first main thm} holds
if the tensor product $M (\bs z)$ is consi\-dered for generic $\bs z$ and
generic (not necessarily diagonal) $Q$. Then the remark in Section~3.1.9
and Theorem~\ref{thm continuity principle}
imply that
the statement of Theorem \ref{The first main thm} holds
for the tensor product of any polynomial
f\/inite-dimensional irreducible $\Y$-modules and any $Q\in GL_N$.
Hence the statement of Theo\-rem~\ref{The first main thm} holds
for direct sums of tensor products of polynomial
f\/inite-dimensional irreducible $\Y$-modules and any $Q\in GL_N$.
\end{proof}

\subsection{Proof of Theorem \ref{thm on completeness}}
\label{prf of thm on completeness}


{\bf 5.6.1.}
For a nonzero weight subspace $M(\bs z)[m_1,\dots,m_N] \subset M(\bs z)$
denote by $d[m_1,\dots,m_N]$ its dimension.
Let $n \geq l_1 \geq \dots \geq l_{N-1} \geq 0$ be the numbers def\/ined in
Section~\ref{Bethe ansatz equations associated with a weight subspace}.

A vector $\bs a = (a_1,\dots,a_n)$
with coordinates $a_i$ from the set $\{0,2,3,\dots,N\}$ will be called
{\it admissible} if for any $j=1,\dots,N-1$ we have \
$l_j\ =\ \# \{\,a_i\ |\ i=1,\dots,n,\,{\rm and}\,\ a_i > j \,\}$
\
In other words, $\bs a$ is admissible if
\ $
m_j\ =\ \# \{\,a_i\ | \ i=1,\dots,n,\ \hbox{and}\ a_i=j \}.
$\

If $\bs a=(a_1,\dots,a_n)$ is admissible, then for $j=1,\dots,N-1$,
there exists a unique increasing map
$\rho_{\bs a,i}\ :\ \{1,\dots,l_i\}\ \to\ \{1,\dots,n\}$
such that $\# \rho_{\bs a,i}^{-1}(j) = 1$ if $a_j > i$
and $\# \rho_{\bs a,i}^{-1}(j) = 0$ if $a_j \leq i$.

We order admissible vectors lexicographically: we say that $\bs a > \bs a'$
if $a_N=a_N',$ $ a_{N-1}=a_{N-1}'$, $\dots$, $a_i=a_i'$,
$a_{i-1} > a_{i-1}$ for some $i$.


\medskip

\noindent
{\bf 5.6.2.}
Let $v = (1,0,\dots,0)\in \C^N$ be the highest weight vector.
Consider the set of vectors
${e_{\bs a}\bs v = e_{a_1,1}v \otimes \dots \otimes e_{a_n,1}v
\in M(\bs z)}$
labeled by admissible indices $\bs a$. Here $e_{a_i,1}v$ denotes $v$
if $a_i = 0$. Then this set of vectors is a basis of
the weight subspace $M(\bs z)[m_1,\dots,m_N]$.
In particular, the total number of admissible indices
equals $d[m_1,\dots,m_N]$.


\medskip

\noindent
{\bf 5.6.3.}
For a function $f(u_1,\dots, u_k)$ set
\[
\Sym_{\,u_1,\dots, u_k}f(u_1,\dots, u_k) =
\sum_{\si\in S_k} f(u_{\si_1},\dots, u_{\si_k}) ,
\]
where the sum is over all permutations $\si$ of $\{1,\dots,k\}$.

\begin{lemma}
The universal weight function $\omega (\bs t,\bs z)$ is given by the rule
\begin{gather}
\label{explicit}
\omega (\bs t,\bs z)=\prod_{s=1}^n\,\prod_{j=1}^{l_1}\big(t_j^{(1)}-z_s\big)
\,\prod_{b=2}^{N-1}\,\prod_{i=1}^{l_{b-1}}\,\prod_{j=1}^{l_b}
\big(t_j^{(b)}-t_i^{(b-1)}\big)
\vspace{1mm}\\
{}\!\times\sum_{\bs a}\,\Sym_{t_1^{(1)},\dots, t_{l_1}^{(1)}}\cdots
\Sym_{t_1^{(N-1)},\dots, t_{l_{N-1}}^{(N-1)}} \Biggl[
\prod_{\satop{s=1}{a_s>1}}^n \Biggl(
\frac1{t_{\rho^{-1}_{\bs a,1}(s)}^{(1)}-z_s}\,\prod_{r=1}^{s-1}\,
\frac{t_{\rho^{-1}_{\bs a,1}(s)}^{(1)}-z_r+1}
{t_{\rho^{-1}_{\bs a,1}(s)}^{(1)}-z_r} \Biggr)
\notag
\vspace{1mm}\\
{}\!\times \prod_{b=2}^{N-1}\,\prod_{\satop{s=1}{a_s>\,b}}^n\Biggl(
\frac1{t_{\rho^{-1}_{\bs a,b}(s)}^{(b)}-t_{\rho^{-1}_{\bs a,b-1}(s)}^{(b-1)}}\,
\prod_{\satop{r=1}{a_r\ge\,b}}^{s-1}
\frac{t_{\rho^{-1}_{\bs a,b}(s)}^{(b)}-t_{\rho^{-1}_{\bs a,b-1}(r)}^{(b-1)}+1}
{t_{\rho^{-1}_{\bs a,b}(s)}^{(b)}-t_{\rho^{-1}_{\bs a,b-1}(r)}^{(b-1)}}
\Biggr)
\,\prod_{c=1}^{N-1}\,\prod_{j=2}^{l_c}\,\prod_{i=1}^{j-1}
\,\frac{t_i^{(c)}-t_j^{(c)}+1}{t_i^{(c)}-t_j^{(c)}}\,
\Biggr]\,e_{\bs a}\bs v ,
\notag
\end{gather}
where the sum is over all admissible $\bs a$\,.
\end{lemma}

The lemma follows from formula~\Ref{cau}, formula~(6.2) in \cite{MTV2},
and Corollaries~3.5, 3.7 in~\cite{TV}.


\medskip

\noindent
{\bf 5.6.4.}
Let $z_1,\dots, z_n$ be real numbers such that $z_{i+1}-z_i>N$. Assume that
$Q=\diag\,(1,q,\dots,\alb q^{N-1})$, where $q$ is a nonzero parameter.
Then the Bethe ansatz equations \Ref{BAE} take the form
\begin{gather}
\label{BAEq}
\prod_{s=1}^n \big(t^{(1)}_j - z_s + 1\big)
\prod_{\satop{j'=1}{j'\neq j}}^{l_1} \big(t^{(1)}_j - t^{(1)}_{j'} - 1 \big)
\prod_{j'=1}^{l_2} \big(t^{(1)}_j - t^{(2)}_{j'}\big)
\\
\qquad{}=
q\,\prod_{s=1}^n\big(t^{(1)}_j - z_s\big) \prod_{\satop{j'=1}{j'\neq j}}^{l_1}
\big(t^{(1)}_j - t^{(1)}_{j'} + 1 \big)
\prod_{j'=1}^{l_2} \big(t^{(1)}_j - t^{(2)}_{j'} - 1\big) ,
\notag
\\
\begin{gathered}
\hspace*{-18mm} \prod_{j'=1}^{l_{a-1}} \big(t^{(a)}_j - t^{(a-1)}_{j'} + 1\big)
\prod_{\satop{j'=1}{j'\neq j}}^{l_a}
\big(t^{(a)}_j - t^{(a)}_{j'} - 1 \big)
\prod_{j'=1}^{l_{a+1}} \big(t^{(a)}_j - t^{(a+1)}_{j'}\big)
\notag
\\
\qquad{}=q\,\prod_{j'=1}^{l_{a-1}} \big(t^{(a)}_j - t^{(a-1)}_{j'}\big)
\prod_{\satop{j'=1}{j'\neq j}}^{l_a}
\big(t^{(a)}_j - t^{(a)}_{j'} + 1 \big)
\prod_{j'=1}^{l_{a+1}} \big(t^{(a)}_j - t^{(a+1)}_{j'}-1\big) ,
\notag
\end{gathered}
\\
\prod_{j'=1}^{l_{N-2}} \big(t^{(a)}_j - t^{(a-1)}_{j'} + 1\big)
\prod_{\satop{j'=1}{j'\neq j}}^{l_{N-1}}
\big(t^{(a)}_j - t^{(a)}_{j'} - 1 \big)
=q\,\prod_{j'=1}^{l_{N-2}} \big(t^{(a)}_j - t^{(a-1)}_{j'}\big)
\prod_{\satop{j'=1}{j'\neq j}}^{l_{N-1}}
\big(t^{(a)}_j - t^{(a)}_{j'} + 1 \big) .
\notag
\end{gather}
For $q=0$ the right hand sides of equations \Ref{BAEq}
equal zero and to solve the equations for $q=0$ one needs to f\/ind
common zeros of the left hand sides.

Let $\bs a = (a_1,\dots,a_N)$ be an admissible index.
Def\/ine $\tilde {\bs t}(\bs a, 0) = (\tilde t^{(i)}_j(\bs a, 0))$
to be the point in $\C^l$ with coordinates
$\tilde t^{(i)}_j(\bs a,0) = z_{\rho_{\bs a,i}(j)}-i$ for all $i$, $j$.
Then $\tilde {\bs t}(\bs a,0)$ is an of\/f-diagonal solution of system
\Ref{BAEq} for $q=0$. That solution has multiplicity one.

Hence for every small nonzero $q$, there exists a unique point
$\tilde {\bs t}(\bs a, q)= (\tilde t^{(i)}_j(\bs a, q))$, such that
\begin{enumerate}\itemsep=0pt
\item[$\bullet$]
$\tilde {\bs t}(\bs a, q)$
is an of\/f-diagonal solution of system \Ref{BAEq} with the same $q$,
\item[$\bullet$]
$\tilde {\bs t}(\bs a, q)$ holomorphically depends on $q$
and tends to $\tilde {\bs t}(\bs a, 0)$ as $q$ tends to zero.
\end{enumerate}
Therefore, for all $i$, $j$, we have
\begin{gather}
\label{easy}
\tilde t^{(i)}_j(\bs a,q) = z_{\rho_{\bs a,i}(j)} - i + O(q) .
\end{gather}

\begin{lemma}
\label{lem on limit of Bethe vectors}
As $q$ tends to zero, the Bethe vector
$\omega(\tilde {\bs t}(\bs a, q),\bs z)$ has the following asymptotics:
\[
\omega(\tilde {\bs t}(\bs a, q),\bs z) =
C_{\bs a}\,e_{\bs a}v + O(q) + \cdots ,
\]
where $C_{\bs a}$ is a nonzero number and the dots denote a linear combination
of basis vectors $e_{\bs a'}v$ with indices $\bs a'$ lexicographically greater
than $\bs a$.
\end{lemma}

The lemma follows from formulae \Ref{easy} and \Ref{explicit}.

\medskip

Lemma \ref{lem on limit of Bethe vectors} implies
Theorem \ref{thm on completeness}.

\section{Comparison theorem}
\label{Sec Comparison theorem}

\subsection[The automorphism $\chi_f$ and the universal difference operator]{The automorphism $\boldsymbol{\chi_f}$ and the universal dif\/ference operator}

\label{auto chi-f}
Let $f(u)$ be a rational function in $u$ whose Laurent expansion at $u=\infty$
has the form $f(u)=1+O(u^{-1})$.
Then the map $\chi_f\,:\,T(u)\, \mapsto\, f(u)\,T(u)$
def\/ines an automorphism of the Yangian $\Y$, see Section~3.1.4. 

For a Yangian module $V$, denote by $V^f$ the representation
of the Yangian on the same vector space given by the rule
\[
X|_{V^f} = (\chi_f(X))|_V
\]
for any $X\in\Y$.

\begin{lemma}
Let $V$ be an irreducible finite-dimensional $\Y$-module with highest weight
series $c_1(u), \dots, c_N(u)$. Then the Yangian module $V^f$ is irreducible
with highest weight series $f(u)c_1(u),\dots,f(u)c_N(u)$.
\end{lemma}

\begin{lemma}
Let $V$ be a finite-dimensional polynomial $\Y$-module.
Then the corresponding transfer matrices satisfy the relation
\[
\Te_a(u)|_{V^f} = f(u)\cdots f(u-a+1) \,\Te_a(u)|_V
\]
for $a = 1,\dots, N$.
\end{lemma}

\subsection{Comparison of kernels}
\label{comparison of kernels}
Let $V$ be a f\/inite-dimensional polynomial $\Y$-module.
Consider the universal dif\/ference operators $\D_{Q,V}(u,\tau)$ and
$\D_{Q,V^f}(u,\tau)$.
Let $C(u)$ be a function satisfying the equation
\[
C(u) = f(u)\, C(u-1) .
\]
Then a $V$-valued function $g(u)$ belongs to the kernel
of $\D_{Q,V}(u,\tau)$ if and only if the function $C(u)g(u)$ belongs to
the kernel of $\D_{Q,V^f}(u,\tau)$.


\medskip

\noindent
{\bf 6.2.1.}
Let $\La = (\La_1,\dots,\La_N)$ be a dominant integral $\glN$-weight and
$a$ a complex number. Denote
$\La(a) = (\La_1 + a , \ldots, \La_N + a)$.
Let $V = M_\La(z)$ be an evaluation module and
\[
f(u) = \frac{u-z}{u-z-a} .
\]
Then $V^f = M_{\La(a)}(z+a)$,
see Section 3.1.8. 

Let $\bs \La = (\La^{(1)},\dots,\La^{(n)})$ be a collection of
integral dominant $\glN$-weights, where $\La^{(i)}=(\La_1^{(i)},\dots,
\La_N^{(i)})$ for $i=1,\dots, n$.
Let $\bs a = (a_1,\dots,a_n)$ be a collection on complex numbers.
Introduce the new collection
$\bs\La(\bs a)=(\La^{(1)}(a_1),\dots,\La^{(n)}(a_n))$.

For complex numbers $z_1,\dots,z_n$, consider
the tensor products of evaluation modules
\begin{gather*}
M_{\bs \La}(\bs z) = M_{\La^{(1)}}(z_1) \otimes \dots
\otimes
M_{\La^{(n)}}(z_n) ,
\\
M_{\bs \La(\bs a)}(\bs z +\bs a) =
M_{\La^{(1)}(a_1)}(z_1+a_1) \otimes\dots\otimes M_{\La^{(n)}(a_n)}(z_n+ a_n) .
\end{gather*}
Let $V = M_{\bs\La}(\bs z)$ and
\[
f(u) = \prod_{s=1}^n \frac{u-z_s}{u-z_s-a_s} .
\]
Then $V^f = M_{\bs\La(\bs a)}(\bs z + \bs a)$.

\begin{theorem}
\label{thm comparison}
An $M_{\bs\La}(\bs z)$-valued function $g(u)$ belongs to the kernel of
$\D_{Q,M_{\bs\La}(\bs z)}(u,\tau)$ if and only if the function $C(u)g(u)$
belongs to the kernel of $\D_{Q,M_{\bs\La(\bs a)}(\bs z + \bs a)}(u,\tau)$,
where
\[
C(u) = \prod_{s=1}^n\, \frac{\Gamma(u-z_s+1)}{\Gamma(u-z_s-a_s+1)} .
\]
\end{theorem}



\noindent
{\bf 6.2.2.}
Assume that $M_{\bs\La}(\bs z)$ is a polynomial Yangian module and
$\bs a = (a_1,\dots,a_N)$ are non-negative integers. Then $C(u)$ is
a polynomial and $M_{\bs\La(\bs a)}(\bs z+\bs a)$ is a polynomial Yangian
module. In this case Theorem \ref{thm comparison} allows us to compare
the quasi-exponential kernels of the dif\/ference operators
$\D_{Q,M_{\bs\La}(\bs z)}(u,\tau)$ and
$\D_{Q,M_{\bs\La(\bs a)}(\bs z + \bs a)}(u,\tau)$.

\section{The kernel in the tensor product of evaluation modules}
\label{The kernel of D}

\subsection{The second main result}
Let $Q=\diag\,(Q_1,\dots,Q_N)\, \in GL_N$.
Let $\bs \La = (\La^{(1)},\dots,$ $\La^{(n)})$ be a collection of
integral dominant $\glN$-weights,
where $\La^{(i)}=(\La_1^{(i)},\dots, \La_N^{(i)})$ and
$\La_N^{(i)}=0$ for $i=1,\dots, n$.

Choose a weight subspace
$M_{\bs \La}(\bs z)[m_1,\dots,m_N]\subset M_{\bs \La}(\bs z)$
and consider the universal dif\/ference operator
$\D_{Q,M_{\bs \La}(\bs z)[\bs m]}$ associated with $Q$ and
$M_{\bs \La}(\bs z)[m_1,\dots,m_N]$.

\begin{theorem}
\label{thm second main, part 1}
Assume that the nonzero numbers $Q_1,\dots,Q_N$ are distinct
and the argument of each of them is chosen.
Then for any $i=1,\dots,N$ and any nonzero vector
$v_0 \in M_{\bs \La}(\bs z)[m_1,\dots,\alb m_N]$, there exists an
$M_{\bs \La}(\bs z)[m_1,\dots,m_N]$-valued quasi-exponential
\begin{gather}
\label{formula of a quasi-exp}
Q_i^u\big(v_0\,u^{m_i} +v_1\,u^{m_i-1} + \dots + v_{m_i}\big)
\end{gather}
which lies in the kernel of
$\D_{Q,M_{\bs \La}(\bs z)[\bs m]}$.
Moreover, all such quasi-exponentials generate the kernel of
$\D_{Q,M_{\bs \La}(\bs z)[\bs m]}$.
\end{theorem}

\begin{proof}
On one hand, by Theorem \ref{The first main thm}, the kernel of
$\D_{Q,M_{\bs \La}(\bs z)[\bs m]}$ is generated by quasi-exponen\-tials.
On the other hand, by Lemma \ref{lem of Q} and part (i) of
Theorem \ref{thm on coefficients}, any quasi-exponential lying in the
kernel of $\D_{Q,M_{\bs \La}(\bs z)[\bs m]}$ must be of the form
\Ref{formula of a quasi-exp}, where
$v_0 $ is a nonzero vector. Since
the kernel is of dimension $N \cdot \dim\, M_{\bs \La}(\bs z)[\bs m]$,
there exists a quasi-exponential of the form
\Ref{formula of a quasi-exp}
with an arbitrary nonzero $v_0$.
\end{proof}



\noindent
{\bf 7.1.1.}
For $i=1,\dots, n$, denote
$S_i = \{ z_i-1,\,z_i - \La^{(i)}_{N-1}-2, \dots,
z_i - \La^{(i)}_{1}-N \}$.

\begin{theorem}
\label{sec local result}
Assume that index $i$ is such that $z_i - z_j \notin \Z$
for any $j\neq i$. Let $f(u)$ be an
$M_{\bs \La}(\bs z)[\bs m]$-valued quasi-exponential
lying in the kernel of
$\D_{Q,M_{\bs \La}(\bs z)[\bs m]}$. Then
\begin{enumerate}\itemsep=0pt
\item[(i)]
$f(u)$
is uniquely determined by its values
\begin{gather*}
v_N=f(z_i-1),\ v_{N-1}=f\big(z_i - \La^{(i)}_{N-1}-2\big), \
\dots,\ v_1 = f\big(z_i - \La^{(i)}_{1} - N\big) \in
M_{\bs \La}(\bs z)[\bs m]
\end{gather*}
at the points of $S_i$.

\item[(ii)]
If $v_N=v_{N-1}=\dots = v_j =0$ for some $j>1$, then $f(z_i-k)=0$ for
$k = 1, 2,\dots,\alb \La^{(i)}_{j-1}+ N-j+1$.

\item[(iii)]
For any vectors $v_N, v_{N-1},\dots, v_1 \in M_{\bs \La}(\bs z)[\bs m]$,
there exists a quasi-exponential $f(u)$ which lies in the kernel of
$\D_{Q,M_{\bs \La}(\bs z)[\bs m]}$ and takes these values at $S_i$.
\end{enumerate}
\end{theorem}

\begin{proof} Consider the polynomial dif\/ference operator{\samepage
\begin{gather*}
\tilde{\D} = \prod_{s=1}^n
\prod_{j=1}^{N-1} (u-z_s-j+1)\,\D_{Q,M_{\bs \La}(\bs z)[\bs m]}
\\
\phantom{\tilde{\D}}{}= \prod_{s=1}^n\prod_{j=1}^{N-1} (u-z_s-j+1) +
\prod_{s=1}^n \prod_{j=2}^{N-1}(u-z_s-j+1)
\,\tilde{ \Te}_1(u)\,\tau^{-1}+\cdots
\\
\phantom{\tilde{\D}=}{}+\prod_{s=1}^n(u-z_s-N+2) \,\tilde{ \Te}_{N-2}(u)\,\tau^{-N+2}
+ \tilde{ \Te}_{N-1}(u)\,\tau^{-N+1}
\\
\phantom{\tilde{\D}=}{}+
(-1)^N\, \prod_{i=1}^N Q_i\, \prod_{s=1}^n
\prod_{i=1}^{N-1}\big(u-z_s+\La^{(s)}_i-i+1\big) \tau^{-N} ,
\end{gather*}
see notation in part (i) of Theorem \ref{thm on coefficients}.
We have $(\tilde \D f)(u) = 0$.}

Condition $(\tilde\D f)(z_i + N - 2)=0$ gives the equation $f(z_i-2)
={\rm const}\, f(z_i-1)$ and thus the value $f(z_i-1)$ determines
the value $f(z_i-2)$. Condition $(\tilde \D f)(z_i+N-3) = 0$ gives
the equation $f(z_i-3) ={\rm const}\, f(z_i-2)$ and hence the
value $f(z_i-2)$ determines the values $f(z_i-3)$. We may continue on
this reasoning up to equation $f(z_i-\La^{(i)}_{N-1}-1) ={\rm
const}\, f(z_i-\La^{(i)}_{N-1})$ which shows that the value
$f(z_i-\La^{(i)}_{N-1})$ determines the value
$f(z_i-\La^{(i)}_{N-1}-1)$.

Condition $(\tilde \D f)(z_i+N-\La^{(i)}_{N-1}-2) = 0$ does not
determine the value $f(z_i-\La^{(i)}_{N-1}-2)$, since
$f(z_i-\La^{(i)}_{N-1}-2)$ enters that condition with coef\/f\/icient $0$.
But condition $(\tilde \D f)(z_i+N-\La^{(i)}_{N-1}-3) {=} 0$ gives the
equation $f(z_i-\La^{(i)}_{N-1}-3) = {\rm const}\,
f(z_i-\La^{(i)}_{N-1}-2)$ which shows that the value
$f(z_i-\La^{(i)}_{N-1}-3)$ is determined by the value
$f(z_i-\La^{(i)}_{N-1}-2)$.

Now we may continue on this reasoning up to the equation
$f(z_i-\La^{(i)}_{N-2}-2) = {\rm const}\, f(z_i-\La^{(i)}_{N-2}-1)$
which determines $f(z_i-\La^{(i)}_{N-2}-2)$ if
$f(z_i-\La^{(i)}_{N-2}-1)$ is known. Condition $(\tilde \D
f)(z_i+N-\La^{(i)}_{N-2}-3) = 0$ does not determine the value
$f(z_i-\La^{(i)}_{N-2}-3)$, but condition $(\tilde \D
f)(z_i+N-\La^{(i)}_{N-2}-4) = 0$ gives the equation
$f(z_i-\La^{(i)}_{N-2}-4) = {\rm const}\,
f(z_i-\La^{(i)}_{N-2}-3)$ and so on. Repeating this reasoning we prove parts
(i) and (ii) of the theorem.

The same reasoning shows that if $f(u) = 0$ for $u\in S_i$, then the
quasi-exponential $f(u)$ identically equals zero. Since the kernel of
$\tilde{\D}$ is generated by quasi-exponentials, we obtain part (iii)
of the theorem.
\end{proof}



\noindent
{\bf 7.1.2.}
Assume that $Q$ is the identity matrix.
Consider the subspace $\sing M_{\bs \La}(\bs z)[\bs m] \!\subset\!
M_{\bs \La}(\bs z)[\bs m]\!$ of $\glN$-singular vectors and the associated
universal dif\/ference operator
$\D_{Q=1,\,\sing M_{\bs \La}(\bs z)[\bs m]}(u,\tau)$.

\begin{theorem}
\label{thm on space for Q=1}
For any $i=1,\dots,N$ and any nonzero vector
$v_0 \in \sing\,M_{\bs \La}(\bs z)[m_1,\dots,m_N]$, there exists
a $\sing\,M_{\bs \La}(\bs z)[m_1,\dots,m_N]$-valued polynomial
\begin{gather}
\label{formula for polynom}
v_0\,u^{m_i + N-i} + v_1\,u^{m_i + N-i-1} +\cdots + v_{m_i+N-i}
\end{gather}
which lies in the kernel of
$\D_{Q=1,\,\sing M_{\bs \La}(\bs z)[\bs m]}$.
Moreover, all such polynomials generate the kernel of
$\D_{Q=1,\, \sing\,M_{\bs \La}(\bs z)[\bs m]}$.
\end{theorem}

\begin{proof}
On one hand, by Theorem \ref{The first main thm}, the kernel of
$\D_{Q=1,\,\sing M_{\bs \La}(\bs z)[\bs m]}$ is generated by
quasi-exponentials. On the other hand, by part (i) of
Theorem \ref{thm on coefficients}, any quasi-exponential lying in the
kernel of $\D_{Q=1,\,\sing M_{\bs \La}(\bs z)[\bs m]}$ must a polynomial.
By Theorem \ref{thm Q=1} such a polynomial has to be of the form indicated in
\Ref{formula for polynom}, where
$v_0 $ is a nonzero vector. Since
the kernel is of dimension $N \cdot \dim\, \sing\,M_{\bs \La}(\bs z)[\bs m]$,
there exists a polynomial of the form \Ref{formula for polynom} with
an arbitrary nonzero~$v_0$.
\end{proof}

\subsection[Kernel of the fundamental difference operator
associated with an eigenvector]{Kernel of the fundamental dif\/ference operator
associated\\ with an eigenvector}
\label{cor on fund spaces}

Assume that $v \in M_{\bs \La}(\bs z)[\bs m]$
is an eigenvector of all transfer matrices,
\[
\Te_{Q,\,M_{\bs \La}(\bs z)[\bs m]}(u)\, v =
\lambda_{k,v}(u) \,v ,
\qquad
k=0,\dots, N .
\]
Then the scalar dif\/ference operator
\[
\Dg_v(u,\tau_u) = \sum_{k=0}^N\,(-1)^k\,\lambda_{k,v}(u)\,\tau_u^{-k}
\]
will be called {\it the fundamental dif\/ference operator}
associated with the eigenvector $v$.

Theorems \ref{thm second main, part 1}, \ref{sec local result}
and~\ref{thm on space for Q=1} give us information on the kernel
of the fundamental dif\/ference operator.

\begin{corollary}\ \ {}
\begin{enumerate}\itemsep=0pt
\item[(i)]
Assume that $Q \in GL_N$ is diagonal with distinct diagonal entries. Then
the kernel of $\D_v(u,\tau_u)$ is generated by
quasi-exponentials $Q_1^up_1(u),\dots, Q_N^up_N(u)$, where for every $i$
the polynomial $p_i(u)\in \C[u]$ is of degree $m_i$.

\item[(ii)] Assume that $Q$ is the identity matrix and the eigenvector
$v$ belongs to
$\sing\,M_{\bs \La}(\bs z)[\bs m]$. Then the kernel of $\D_v(u,\tau_u)$
is generated by suitable polynomials $p_1(u),\dots,p_N(u)$ of degree
$m_1+N-1,\dots, m_N$, respectively.
\end{enumerate}
\end{corollary}

\begin{corollary}
Assume the index $i$ is such that $z_i - z_j \notin \Z$ for any $j\neq i$.
Let $f(u)$ be a~quasi-exponential lying in the kernel of $\D_{v}(u,\tau_u)$.
Then
\begin{enumerate}\itemsep=0pt
\item[(i)]
$f(u)$
is uniquely determined by its values
\[
v_N=f(z_i-1),\quad v_{N-1}=f(z_i - \La^{(i)}_{N-1}-2),
\quad \dots,\quad v_1 = f(z_i - \La^{(i)}_{1} - N)
\]
at the points of $S_i$.

\item[(ii)]
If $v_N=v_{N-1}=\dots = v_j =0$ for some $j>1$, then $f(z_i - k)=0$ for
$k = 1, 2,\dots, \La^{(i)}_{j-1}+ N-j+1$.

\item[(iii)]
For any numbers $v_N,v_{N-1},\dots,v_1$ there exists a quasi-exponential $f(u)$
which lies in the kernel of $\D_{v}(u,\tau_u)$ and takes these values at $S_i$.
\end{enumerate}
\end{corollary}

If $v$ is a Bethe eigenvector, then these two corollaries were proved in
\cite{MV2} and \cite{MV3}.

\section{Generating operator of the Gaudin transfer matrices}
\label{Sec Gen Gaudin}

\subsection[Current algebra $\glN\lbrack x\rbrack$]{Current algebra $\boldsymbol{\glN[x]}$}


{\bf 8.1.1.}
Let $\glN[x]$ be the Lie algebra of polynomials with coef\/f\/icients in $\glN$
with point-wise commutator. The elements $e^{\{s\}}_{ab} = e_{ab}x^s$ with
$a,b = 1,\dots,N$, $s=0,1,\dots,$ span $\glx$. We have
$[e^{\{r\}}_{ab},e^{\{s\}}_{cd}] = \delta_{bc} e^{\{r+s\}}_{ad}
- \delta_{ad} e^{\{r+s\}}_{cb}$. We shall identify $\glN$ with
the subalgebra of $\glN[x]$ of constant polynomials.

For $a,b = 1,\dots, N$, we set
\[
L_{ab}(u) = \sum_{s=0}^\infty e^{\{s\}}_{ba} u^{-s-1}
\]
and $L(u) = \sum_{a,b} E_{ab} \otimes L_{ab}(u) \in \End\,(\C^N) \otimes
\glN[x]((u^{-1}))$.

Let $V$ be a $\glN$-module. For $z\in \C$ denote by $M(z)$ the corresponding
evaluation $\glN[x]$-module on the same vector space, where we def\/ine
$e^{\{s\}}_{ab}|_V = z^se_{ab}$ for all $s, a,b$. Then the
series $L_{ab}(u)|_V$ converges to the $\End\,(V)$-valued
function $(u-z)^{-1} e_{ba}$,
which is a rational function in~$u$.

Let $\La=(\La_1,\dots,\La_N)$ be an integral dominant $\glN$-weight,
$M_\La$ the corresponding irreducible highest weight $\glN$-module,
$M_\La(z)$ the associated evaluation $\glN[x]$-module. Then the
module $M_\La(z)$ is called {\it polynomial} if $\La_N$ is a non-negative
integer.


\medskip

\noindent
{\bf 8.1.2.}
Let $\bs \La = (\La^{(1)},\dots,\La^{(n)})$ be a collection of
integral dominant $\glN$-weights,
where
$\La^{(i)}=(\La_1^{(i)},\dots, \La_N^{(i)})$ for $i=1,\dots, n$.
For $z_1,\dots,z_n\,\in\,\C$, we consider
the tensor product of $\glN[x]$-evaluation modules:
\[
\M_{\bs \La}(\bs z)
= M_{\La^{(1)}}(z_1) \otimes \dots
\otimes
M_{\La^{(n)}}(z_n) .
\]
For any $a,b$, the series $L_{ab}(u)$ acts on $\M_{\bs \La}(\bs z)$
by the formula
\[
L_{ab}(u)|_{\M_{\bs\La}(\bs z)} = \sum_{j=1}^n\ \frac{e^{(j)}_{ba}}{u-z_j} .
\]
If $z_1,\dots,z_N$ are distinct, then $\M_{\bs \La}(\bs z)$
is an irreducible $\glN[x]$-module.
Let
\[
M_{\La^{(1)}} \otimes \dots\otimes M_{\La^{(n)}} =
\oplus_\La\, M_\La
\]
be the decomposition of the tensor product of $\glN$-modules into the direct
sum of irreducible $\glN$-modules. Then for any $z\in \C$,
\[
M_{\La^{(1)}}(z) \otimes \dots\otimes M_{\La^{(n)}}(z) =
\oplus_\La\, M_\La (z)
\]
is the decomposition of the tensor product of evaluation $\glN[x]$-modules
into the direct sum of irreducible $\glN[x]$-modules.


\medskip

\noindent
{\bf 8.1.3.}
Consider $\C^N$ as the $\glN$-module with highest weight
$(1,0,\dots,0)$. Let $\La =$
$(\La_1,\dots,\La_N)$ be an integral dominant $\glN$-weight
with $\La_N \in \Z_{\geq 0}$. Then there exists $k\in \Z_{\geq 0}$ such
that~$(\C^N)^{\otimes k}$ contains $M_\La$ as a $\glN$-submodule.

The previous remarks show that for
any $z_1,\dots,z_n \in\C$, all irreducible submodules of
the $\glN[x]$-module $\C^N(z_1)\otimes \dots \otimes \C^N(z_n)$ are
tensor products of polynomial evaluation $\glN[x]$-modules.
Moreover, if $V$ is a $\glN[x]$-module which is the tensor product
of polynomial evaluation
$\glN[x]$-modules, then there exist $z_1,\dots,z_n$ such that
$V$ is isomorphic to a submodule of the $\glN[x]$-module
$\C^N(z_1)\otimes \dots \otimes \C^N(z_n)$.


\medskip

\noindent
{\bf 8.1.4.}
Let ${\pi:U(\glN)\to\End\,(\C^N)}$ be the representation homomorphism for
the $\glN$-mo\-du\-le~$\C^N$. Clearly, for any $x\in U(\glN)$ we have
\begin{gather}
\label{adjoint Gaudin}
\bigl[\pi(x)\otimes 1+1\otimes x ,\,L(u)\bigr] = 0.
\end{gather}

For a non-degenerate matrix $A\in\End\,(\C^N)$,
def\/ine an automorphism $\nu_A$ of $\glN[x]$ by the formula
\[
(\id\otimes \nu_A)\bigl(L(u)\bigr) =
\sum_{ab}\ A^{-1}\!E_{ab}\>A\otimes L_{ab}(u) .
\]

Let $V$ be a f\/inite-dimensional $\glN[x]$-module with the representation
$\mu:\glN[x] \to\End\,(V)$ and $\tilde\mu:GL_N\to\End\,(V)$ the corresponding
representation of the group $GL_N$. The automorphism $\nu_A$ induces
a new $\glN[x]$-module structure $V^A$ on the same vector space with
the representation $\mu_A=\mu\circ\nu_A$. Formula \Ref{adjoint Gaudin}
yields that for any $x\in\glN[x]$,
\[
\mu_A(x) = \tilde\mu(A)\,\mu(x)\bigl(\tilde\mu(A)\bigr)^{-1} ,
\]
that is, the $\glN[x]$-modules $V$ and $V^A$ are isomorphic. In particular,
if $V$ is the tensor product of polynomial evaluation $\glN[x]$-modules, then
$V^A$ is the
tensor product of polynomial evaluation $\glx$-modules too.

\subsection[Fundamental differential operator]{Fundamental dif\/ferential operator}

Let $K=(K_{ab})$ be an $N\times N$ matrix with complex entries.
For $a,b=1,\dots,N$, def\/ine the dif\/ferential operator{\samepage
\[
X_{ab}(u,\p) = \delta_{ab}\, \p - K_{ab} - L_{ab}(u) ,
\]
where $\p=d/du$.}

Following \cite{T}, introduce the dif\/ferential operator
\[
\Dg(u,\p) = \sum_{\si\in S_N} (-1)^\si\,
X_{1\,\si_1}(u,\p)\,X_{2\,\si_2}(u,\p)\cdots X_{N\si_N}(u,\p) ,
\]
where the sum is over all permutations $\si$ of $\{1,\dots,N\}$.
The operator $\Dg(u,\p)$ will be called {\it the universal differential
operator\/} associated with the matrix $K$.

\begin{lemma}
\label{rowdet2}
Let $\pi$ be a map $\{1,\dots,N\}\to\{1,\dots,N\}$.
If $\pi$ is a permutation of $\{1,\dots,N\}$, then
\[
\sum_{\si\in S_N} (-1)^\si\,
X_{\pi_1\si_1}(u,\tau)\,X_{\pi_2\si_2}(u,\tau)\cdots
X_{\pi_N\si_N}(u,\tau) = (-1)^\pi\,\Dg(u,\p) .
\]
If $\pi$ is not bijective, then
\[
\sum_{\si\in S_N} (-1)^\si\,
X_{\pi_1\si_1}(u,\tau)\,X_{\pi_2\si_2}(u,\tau)\cdots
X_{\pi_N\si_N}(u,\tau) = 0 .
\]
\end{lemma}

The statement is Proposition~8.1 in~\cite{MTV2}.


\medskip

\noindent
{\bf 8.2.1.}
Introduce the coef\/f\/icients $\Tee_0(u),\dots, \Tee_N(u)$ of
$\Dg(u,\p)$:
\[
\Dg (u,\p) = \sum_{k=0}^N\,(-1)^k\,\Tee_{k}(u)\,\p^{N-k} ,
\]
in particular, $\Tee_0(u) = 1$.
The coef\/f\/icients $\Tee_k(u)$ are called {\it the transfer matrices of
the Gaudin type model} associated with $K$.

The main properties of the transfer matrices:
\begin{enumerate}\itemsep=0pt
\item[(i)]
the transfer matrices commute:
$
[\Tee_k(u),\Tee_l(v)] = 0
$
for all $k$, $l$, $u$, $v$,
\item[(ii)]
if $K$ is a diagonal matrix, then the transfer matrices preserve
the $\glN$-weight:
$
[\Tee_k(u), e_{aa}] {=} 0
$
for all $k$, $a$, $u$,
\item[(iii)]
if $K$ is the zero matrix, then the transfer matrices commute with
the subalgebra $U(\glN)$:
$
[\Tee_k(u), x] = 0
$
for all $k,u$ and $x \in U(\glN)$,
\end{enumerate}
see \cite{T, MTV2}.


\medskip

\noindent
{\bf 8.2.2.}
If $V$ is the tensor product of evaluation f\/inite-dimensional
$\glN[x]$-modules, then the universal operator $\Dg(u,\p)$
induces a dif\/ferential operator acting on $V$-valued functions in
$u$. This operator
will be called {\it the universal differential
operator} associated with $K$ and $V$ and denoted by
$\Dg_{K,V}(u,\p)$. The linear
operators $\Tee_k(u)|_{V} \in \End\,(V)$ will be
called {\it the transfer matrices} associated with $K$ and
$V$ and denoted by $\Tee_{k,K,V}(u)$.
They are rational functions in $u$.


\medskip

\noindent
{\bf 8.2.3.}
If $\Dg_K(u,\tau)$ is the universal dif\/ferential operator associated with
the matrix $K$ and $\nu_A : \glx \to \glx$ is the automorphism def\/ined in
Section 8.1.4. 
Then Lemma~\ref{rowdet2} implies that
\[
\nu_A(\Dg_K(u,\p)) = \Dg_{AKA^{-1}}(u,\p)
\]
is the universal dif\/ferential operator associated with the matrix $AKA^{-1}$,
cf.~Lemma \ref{sec on nu-A 2}.

Let $V$ be the tensor product of f\/inite-dimensional evaluation
$\glN[x]$-modules, and $\tilde\mu:G_N\to GL(V)$ the associated
$GL_N$-representation. Then
\[
\Dg_{AKA^{-1}}(u,\p)|_V =
\tilde\mu(A)\,\Dg_K(u,\p)|_{V}\,\tilde\mu \big(A^{-1}\big) .
\]

\subsection{More properties of the Gaudin type transfer matrices}

Let $K = \diag\,(K_1, \dots, K_N)$.
Let $\bs \La = (\La^{(1)},\dots,\La^{(n)})$ be a collection of
integral dominant $\glN$-weights,
where $\La^{(i)}=(\La_1^{(i)},\dots, \La_N^{(i)})$ and
$\La_N^{(i)}=0$ for $i=1,\dots, n$.

For $\bs m = (m_1,\dots,m_N)$ denote by
$\M_{\bs \La}(\bs z)[\bs m] \,\subset\, \M_{\bs \La}(\bs z)$ the
weight subspace of the $\glN$-weight $\bs m$ and by
$\sing\, \M_{\bs \La}(\bs z)[\bs m] \,
\subset\, \M_{\bs \La}(\bs z)[\bs m]$ the subspace of $\glN$-singular vectors.

Consider the universal dif\/ferential operator
\[
\Dg_{K,\,\M_{\bs \La}(\bs z)}
(u,\p) = \sum_{k=0}^N\,(-1)^k\,
\Tee_{k, K,\,\M_{\bs \La}(\bs z)}(u)
\,\p^{N-k}
\]
associated with $K$ and $\M_{\bs \La}(\bs z)$.
We have $\Tee_{0,K,\,\M_{\bs \La}(\bs z)[\bs m]}(u)= 1$ and
\[
\Tee_{k,K,\,\M_{\bs \La}(\bs z)[\bs m]}(u) = \Tee_{k0} + \Tee_{k1}u^{-1}
+ \Tee_{k2}u^{-2} + \cdots
\]
for suitable $\Tee_{ki} \in \End\,(\M_{\bs \La}(\bs z)[\bs m])$.

\begin{theorem}
\label{thm on coefficients Gaudin}
The following statements hold.
\begin{enumerate}\itemsep=0pt
\item[(i)]
The operators $\Tee_{10},\Tee_{20},\dots,\Tee_{N0}$ and
$\Tee_{11},\Tee_{21},\dots,\Tee_{N1}$ are scalar operators.
Moreover, the following relations hold:
\begin{gather*}
x^N+\,\sum_{k=1}^N \,(-1)^{k}\,\Tee_{k0}\,x^{N-k} = \prod_{i=1}^N\,(x-K_i),
\\
\sum_{k=1}^N \,(-1)^{k}\,\Tee_{k1}\,x^{N-k} =
- \prod_{i=1}^N\,(x-K_i) \sum_{j=1}^N\,\frac{m_j}{x-K_j} .
\end{gather*}

\item[(ii)]
For $k=1,\dots,N-1$, we have
\begin{gather}
\label{formula on Polyn}
\Tee_{k,Q,\M_{\bs \La}(\bs z)}(u) =
\tilde{\Tee}_{k}(u) \prod_{s=1}^n\frac 1{(u-z_s)^k} ,
\end{gather}
where $\tilde{\Tee}_{k}(u)$ is a polynomial in $u$ of degree $nk$. Moreover,
the operators
\[
\bar\Tee_{k,r} = \tilde{\Tee}_{k}(z_r)\!\!
\prod_{s=1, s\ne r}^n\frac1{(z_r-z_s)^k}
\]
are scalar and
\begin{gather}
\label{chara eqn}
\sum_{k=0}^N \,(-1)^{k} \,\bar\Tee_{k,s}\!\! \prod_{j=0}^{N-r-1} (d-j) =
\prod_{i=1}^N \big(d-\La^{(r)}_i-N+i\big) .
\end{gather}
\end{enumerate}
\end{theorem}

\begin{proof}
Part (i) follows from Proposition B.1 in \cite{MTV2}. The existence of
presentation \Ref{formula on Polyn} follows from the def\/inition of the
universal dif\/ferential operator. To prove equation \Ref{chara eqn} it is enough
to notice that the leading singular term of the universal dif\/ferential operator
at $u=z_r$ is equal to the leading singular term of the universal dif\/ferential
operator associated with one evaluation module $M_{\La^{(r)}}(z_r)$, which in
its turn expresses via the quantum determinant.
\end{proof}


\medskip

\noindent
{\bf 8.3.1.}
Assume that $K$ is the zero matrix. Then the associated transfer matrices
preserve $\sing\,\M_{\bs \La}(\bs z)[\bs m]$ and we may consider the universal
dif\/ferential operator
\[
\Dg_{K=0,\,\sing\M_{\bs \La}(\bs z)[\bs m]}(u,\p)
= \sum_{k=0}^N \, (-1)^k\, {\Tee}_{k}(u) \,\p^{N-k}
\]
acting on $\sing\,\M_{\bs \La}(\bs z)[\bs m]$-valued functions of $u$.

\begin{theorem}
\label{thm K=0}
For $k=1,\dots, N$, the coefficients $ {\Tee}_{k}(u)$ have the following
Laurent expansion at $u=\infty$:
\[
{\Tee}_{k}(u) = {\Tee}_{k,0} u^{-k} + {\Tee}_{k,1} u^{-k-1} + \cdots ,
\]
where the operators ${\Tee}_{1,0},\dots,{\Tee}_{N,0}$ are scalar operators.
Moreover,
\[
\sum_{k=0}^N (-1)^k \Tee_{k,0} \!\! \prod_{j=0}^{N-k-1} (d-j)
= \prod_{s=1}^N\,(d - m_s - N + s) .
\]
\end{theorem}

The proof of Theorem~\ref{thm K=0} is similar to the proof of
Theorem~\ref{thm Q=1}.

\subsection{First main result in the Gaudin case}

\begin{theorem}
\label{The first main thm Gaudin}
Let $K$ be an $N\times N$-matrix and $\M_{\bs\La}(\bs z)$ the tensor product
of polynomial $\glx$-modules. Consider the universal differential operator
$\Dg_{K,\,\M_{\bs\La}(\bs z)}(u,\p)$ associated with~$K$ and~$\M_{\bs\La}(\bs z)$. Then the kernel of
$\Dg_{K,\,\M_{\bs\La}(\bs z)}(u,\p)$ is generated by quasi-exponentials.
\end{theorem}

A statement of this type was conjectured in \cite{CT}.
Theorem \ref{The first main thm Gaudin} will be proved
in Section \ref{Bethe ansatz Gaudin}.

\section[Continuity principle for differential
operators with quasi-exponential kernel]{Continuity principle for dif\/ferential
operators\\ with quasi-exponential kernel}
\label{Sec Contin Gaudin}

\subsection{Quasi-exponentials}
\label{sec Quasi-exponentials Gaudin}

Let $V$ be a complex vector space of dimension $d$. Let
$A_0(u),\dots, A_N(u)$ be $\End\,(V)$-valued rational functions
in $u$. Assume that each of these functions has limit as $u \to \infty$
and $A_0(u) = 1$ in $\End\,(V)$. Then the dif\/ferential operator
\[
\Dg = \sum_{k=0}^N\,A_k(u)\, \p^{N-k} ,
\]
acting on $V$-valued functions in $u$, will be called admissible at inf\/inity.

For every $k$, let $A_k(u) = A^\infty_{k,0} + A^\infty_{k,1}
u^{-1} + A^\infty_{k,2} u^{-2} + \cdots{}$
be the Laurent expansion at inf\/inity. Consider the algebraic equation
\begin{gather}
\label{char eqn Gaudin}
\det\, \left(x^N + x^{N-1} A_{1,0} + \dots +
x A_{N-1,0} + A_{N,0} \right) = 0
\end{gather}
with respect to variable $x$.

\begin{lemma}
\label{lem of Q Gaudin}
If a nonzero $V$-valued quasi-exponential
$e^{\la u}(u^dv_d + u^{d-1}v_{d-1} + \dots + v_0)$ lies in the kernel
of an admissible at infinity differential operator $\Dg$, then
$\la$ is a root of equation \Ref{char eqn Gaudin}.
\end{lemma}

\subsection{Continuity principle}
Let $A_0(u,\epsilon),\dots, A_N(u,\ep)$ be $\End\,(V)$-valued rational
functions in $u$ analytically depending on $\ep \in [0,1)$. Assume that
\begin{enumerate}\itemsep=0pt
\item[$\bullet$]
for every $\ep \in [0,1)$ the dif\/ference operator
$\Dg_\ep = \sum\limits_{k=0}^N A_k(u,\ep)\, \p^{N-k}$
is admissible at inf\/inity,

\item[$\bullet$]
for every $\ep \in (0,1)$ the kernel of $\Dg_\ep$ is
generated by quasi-exponentials,
\item[$\bullet$]
there exists a natural number $m$ such that for every $\ep \in (0,1)$ all
quasi-exponentials generating the kernel of $\Dg_\ep$ are of degree less than~$m$.
\end{enumerate}

\begin{theorem}
\label{thm continuity principle Gaudin}
Under these conditions the kernel of
the differential operator $\Dg_{\ep=0}$ is generated
by quasi-exponentials.
\end{theorem}

The proof is similar to the proof of Theorem~\ref{thm continuity principle}.

\section {Bethe ansatz in the Gaudin case}
\label{Bethe ansatz Gaudin}

\subsection{Preliminaries}

Consider $\C^N$ as the $\glN$-module with highest weight
$(1,0,\dots,0)$. For complex numbers $z_1,\dots,\alb z_n$, denote
\[
\M(\bs z) = \C^N(z_1)\otimes \dots \otimes \C^N(z_n) ,
\]
which is the tensor product of polynomial $\glx$-modules.
Let
\[
\M(\bs z) = \oplus_{m_1\geq\dots \geq m_N\geq 0}\, \M(\bs z)[\bs m]
\]
be its $\glN$-weight decomposition with respect to the Cartan subalgebra of
diagonal matrices. Here $\bs m = (m_1,\dots,m_N)$.

Assume that $K = \diag\,(K_1,\dots,K_N)$ is a diagonal
$N\times N$-matrix with distinct coordinates
and consider the universal dif\/ferential operator
\[
\Dg_{K,\,\M(\bs z)}(u,\p)=
\sum_{k=0}^N \,(-1)^k\,\Tee_{k,K,\, \M(\bs z)}(u)\,\p^{N-k}.
\]
associated with $\M(\bs z)$ and $K$.
Acting on $\M(\bs z)$-valued functions, the operator
$\Dg_{K,\,\M(\bs z)}(u,\p)$ preserves the weight decomposition.

In this section we shall study the kernel of this operator restricted
to $\M(\bs z)[\bs m]$-valued functions.

\subsection{Bethe ansatz equations associated with a weight subspace}
\label{Bethe ansatz equations associated with a weight subspace Gaudin}
Consider a nonzero weight subspace $\M(\bs z)[\bs m]$. Introduce
$\bs l = (l_1,\dots,l_{N-1})$ with
$
l_j = m_{j+1} + \dots + m_N.
$
We have $n\geq l_1\geq \dots \geq l_{N-1}\geq 0$.
Set $l_0=l_N=0$ and $l=l_1 + \dots + l_{N-1}$.
In what follows we shall consider functions of $l$ variables
\[
\bs t =
\big(t^{(1)}_{1},\dots, t^{(1)}_{l_1},
t^{(2)}_{1},\dots, t^{(2)}_{l_2},\dots,
t^{(N-1)}_{1},\dots, t^{(N-1)}_{l_{N-1}} \big) .
\]
The following system of $l$ algebraic equations with respect to $l$
variables $\bs t$ is called {\it the Bethe ansatz equations}
associated with $\M(\bs z)[\bs m]$ and matrix $K$,
\begin{gather}
\label{BAE Gaudin}
\sum_{s=1}^n
\frac 1
{t^{(1)}_j - z_s}
+
\sum_{j'=1}^{l_2}
\frac 1{t^{(1)}_j - t^{(2)}_{j'}}
-
\sum_{\satop{j'=1}{j'\neq j}}^{l_1}
\frac 2
{t^{(1)}_j - t^{(1)}_{j'}}=
K_2-K_1 ,
\\
\sum_{j'=1}^{l_{a+1}}
\frac 1{t^{(a)}_j - t^{(a+1)}_{j'}}
+
\sum_{j'=1}^{l_{a-1}}
\frac 1{t^{(a)}_j - t^{(a-1)}_{j'}}
-
\sum_{\satop{j'=1}{j'\neq j}}^{l_a}
\frac 2
{t^{(a)}_j - t^{(a)}_{j'}}=
K_{a+1}-K_{a} ,\nonumber\\
\sum_{j'=1}^{l_{N-2}}
\frac 1{t^{(N-1)}_j - t^{(N-2)}_{j'}}
-
\sum_{\satop{j'=1}{j'\neq j}}^{l_{N-1}}
\frac 2
{t^{(N-1)}_j - t^{(N-1)}_{j'}}=
K_{N}-K_{N-1} .\nonumber
\end{gather}
Here the equations of the f\/irst group are labeled by $j=1,\dots,l_1$,
the equations of the second group are labeled by $a=2,\dots,N-2$,
$j=1,\dots,l_a$, the equations of
the third group are labeled by $j=1,\dots,l_{N-1}$.

\subsection{Weight function and Bethe ansatz theorem}

We denote by $\omega (\bs t,\bs z)$ {\it the universal weight function
of the Gaudin type} associated with the weight subspace $\M(\bs z)[\bs
m]$. The universal weight function of the Gaudin type is def\/ined in
\cite{SV}. A convenient formula for $\omega (\bs t,\bs z)$ is given
in Appendix in \cite{RSV} and Theorems 6.3 and 6.5 in \cite{RSV}.

If $\tilde{\bs t}$ is a solution of the Bethe ansatz equations, then
the vector $\omega(\tilde {\bs t},\bs z) \in \M(\bs z)[\bs m]$
is called {\it the Bethe vector} associated with $\tilde{\bs t}$.

\begin{theorem}
\label{thm on Bethe ansatz Gaudin}
Let $K$ be a diagonal matrix and $\tilde{\bs t}$ a solution
of the Bethe ansatz equa\-tions \Ref{BAE Gaudin}.
Assume that the Bethe vector $\omega(\tilde {\bs t},\bs z)$ is nonzero.
Then the Bethe vector is an eigenvector of all transfer-matrices
$\Tee_{k, K, \M(\bs z)}(u)$, $k=0,\dots,N$.
\end{theorem}

The statement follows from Theorem~9.2 in~\cite{MTV2}.
For ${k=1}$, the result is established in~\cite{RV}.

\medskip

The eigenvalues of the Bethe vector are given by the following construction.
Set
\begin{gather*}
\chi^1(u,\bs t,\bs z)
=
K_1 + \sum_{s=1}^n \frac 1 {u-z_s}
- \sum_{i=1}^{l_1}
\frac 1 {u- t^{(1)}_i} ,
\\
\chi^a(u, \bs t,\bs z) =
K_a +
\sum_{i=1}^{l_{a-1}}
\frac 1 {u- t^{(a-1)}_i}
- \sum_{i=1}^{l_a}
\frac 1 {u- t^{(a)}_i} ,
\end{gather*}
for $ a=2,\dots, N$. Def\/ine the functions
$\lambda_k(u,\bs t,\bs z)$ by the formula
\[
\big(\p -\chi^1(u,\bs t,\bs z)\big)\cdots
\big(\p-\chi^N(u,\bs t,\bs z)\big) =
\sum_{k=0}^N \,(-1)^k\,\lambda_k(u,\bs t,\bs z)\,\p^{N-k} .
\]
Then for $k=0,\dots, N$,
\[
\Tee_{k, Q, M(\bs z)}(u)\,\omega(\tilde{\bs t},\bs z) =
\lambda_k(u,\tilde{\bs t},\bs z)\,
\omega(\tilde{\bs t},\bs z) ,
\]
see Theorem~9.2 in \cite{MTV2}.

\subsection[Differential operator associated with a solution of
the Bethe ansatz equations]{Dif\/ferential operator associated with a solution\\ of
the Bethe ansatz equations}

Let $\tilde{\bs t}$ be a solution of system \Ref{BAE Gaudin}.
The scalar dif\/ference operator
\[
\Dg_{\tilde {\bs t}}(u,\tau) = \sum_{k=0}^N \,(-1)^k\,
\lambda_k(u,\tilde{\bs t},\bs z)\,\p^{N-k}
\]
will be called the associated {\it fundamental differential operator}.

\begin{theorem}
\label{thm on fund operator Gaudin}
The kernel of $\Dg_{\tilde {\bs t}}(u,\tau)$ is generated by quasi-exponentials
of degree bounded from above by a function in $n$ and $N$.
\end{theorem}

This is Proposition~6.4 in~\cite{MV3},
which is a generalization of Lemma~5.6 in~\cite{MV1}.

\subsection{Completeness of the Bethe ansatz}
\begin{theorem}
\label{thm on completeness Gaudin}
Let $z_1,\dots,z_n$ and $K=\diag\,(K_1,\dots,K_N)$ be generic.
Then the Bethe vectors form a basis in $\M(\bs z)[\bs m]$.
\end{theorem}

Theorem \ref{thm on completeness Gaudin} will be proved
in Section \ref{prf of thm on completeness Gaudin}.

\begin{corollary}
Theorems {\rm \ref{thm on completeness Gaudin}}
and~{\rm \ref{thm continuity principle Gaudin}} imply
Theorem {\rm \ref{The first main thm Gaudin}}.
\end{corollary}

\begin{proof}
Theorems \ref{thm on fund operator Gaudin}
and \ref{thm on completeness Gaudin} imply that
the statement of Theorem \ref{The first main thm Gaudin} holds
if the tensor product $\M (\bs z)$ is considered for generic $\bs z$ and
generic diagonal $K$. Then
according to the remark in Section 8.2.3, 
the statement of Theorem \ref{The first main thm} holds
if the tensor product $\M (\bs z)$ is considered for generic $\bs z$ and
generic (not necessarily diagonal) $K$. Then the remark in Section~8.1.3
and Theorem \ref{thm continuity principle Gaudin} imply that
the statement of Theorem \ref{The first main thm Gaudin} holds
for the tensor product of any polynomial f\/inite-dimensional
$\glx$-modules and any $K$.
\end{proof}

\subsection{Proof of Theorem \ref{thm on completeness Gaudin}}
\label{prf of thm on completeness Gaudin}


{\bf 10.6.1.}
For a nonzero weight subspace $\M(\bs z)[\bs m] \subset \M(\bs z)$
denote by $d[\bs m]$ its dimension.
Let $n \geq l_1 \geq \dots \geq l_{N-1} \geq 0$ be the numbers def\/ined in
Section~\ref{Bethe ansatz equations associated with a weight subspace Gaudin}.

A vector $\bs a = (a_1,\dots,a_n)$
with coordinates $a_i$ from the set $\{0,2,3,\dots,N\}$ will be called
{\it admissible} if for any $j=1,\dots,N-1$ we have
$l_j = \# \{\,a_i\ |\ i=1,\dots,n,\,{\rm and}\,\ a_i > j \,\}$.
In other words, $\bs a$ is admissible if
$m_j = \# \{\,a_i\ | \ i=1,\dots,n,\ \hbox{and}\ a_i=j \}$.

If $\bs a = (a_1,\dots,a_n)$ is admissible, then for $j=1,\dots,N-1$,
there exists a unique increasing map
$\rho_{\bs a,i}\ :\ \{1,\dots,l_i\}\ \to\ \{1,\dots,n\}$
such that $\# \rho_{\bs a,i}^{-1}(j) = 1$ if $a_j > i$
and $\# \rho_{\bs a,i}^{-1}(j) = 0$ if $a_j \leq i$.

Let $v = (1,0,\dots,0)\in \C^N$ be the highest weight vector.
The set of vectors
$
e_{\bs a}\bs v = e_{a_1,1}v \otimes \dots \otimes e_{a_n,1}v
\in \M(\bs z)$,
labeled by admissible indices $\bs a$, form a basis of
$\M(\bs z)[\bs m]$.


\medskip

\noindent
{\bf 10.6.2.}
Let $z_1,\dots, z_n$ be distinct numbers. Assume that
$K = \diag\,(1/q,2/q,\dots,N/q)$, where $q$ is a small nonzero parameter.
Then the Bethe ansatz equations \Ref{BAE Gaudin} take the form
\begin{gather}
\label{BAEq Gaudin}
\sum_{s=1}^n
\frac 1
{t^{(1)}_j - z_s}
+
\sum_{j'=1}^{l_2}
\frac 1{t^{(1)}_j - t^{(2)}_{j'}}
-
\sum_{\satop{j'=1}{j'\neq j}}^{l_1}
\frac 2
{t^{(1)}_j - t^{(1)}_{j'}} =
\frac 1q ,
\\
\sum_{j'=1}^{l_{a+1}}
\frac 1{t^{(a)}_j - t^{(a+1)}_{j'}}
+
\sum_{j'=1}^{l_{a-1}}
\frac 1{t^{(a)}_j - t^{(a-1)}_{j'}}
-
\sum_{\satop{j'=1}{j'\neq j}}^{l_a}
\frac 2
{t^{(a)}_j - t^{(a)}_{j'}} =
\frac 1q ,\nonumber \\
\sum_{j'=1}^{l_{N-2}}
\frac 1{t^{(N-1)}_j - t^{(N-2)}_{j'}}
-
\sum_{\satop{j'=1}{j'\neq j}}^{l_{N-1}}
\frac 2
{t^{(N-1)}_j - t^{(N-1)}_{j'}}=
\frac 1q .\nonumber
\end{gather}

\begin{lemma}
For any admissible index $\bs a = (a_1,\dots,a_N)$ and small nonzero $q$,
there exists a~solution $\tilde {\bs t}(\bs a, q)$
$= (\tilde t^{(i)}_j(\bs a, q))$ of system \Ref{BAEq Gaudin}
such that
\begin{gather}
\label{easy Gaudin}
\tilde t^{(i)}_j(\bs a,q) = z_{\rho_{\bs a,i}(j)} +
\sum_{k=1}^{i}\, \frac q{a_{\rho_{\bs a,i}(j)}-k}
+ O\big( q^2\big)
\end{gather}
for every $i$, $j$.
\end{lemma}

\begin{lemma}
\label{lem on limit of Bethe vectors Gaudin}
As $q$ tends to zero, the Bethe vector
$\omega(\tilde {\bs t}(\bs a, q),\bs z)$ has the following asymptotics:
\[
\omega(\tilde {\bs t}(\bs a, q),\bs z) =
C_{\bs a}\,q^{-l} e_{\bs a}v + O\big(q^{-l+1}\big) ,
\]
where $C_{\bs a}$ is a nonzero number.
\end{lemma}

The lemma follows from formula \Ref{easy Gaudin} and the formula
for the universal weight function in~\cite{RSV}.

\medskip

Lemma \ref{lem on limit of Bethe vectors Gaudin} implies Theorem
\ref{thm on completeness Gaudin}.

\section{Comparison theorem in the Gaudin case}
\label{Comparison theorem in the Gaudin case}

\subsection[The automorphism $\tilde\chi_f$ and
the universal differential operator]{The automorphism $\boldsymbol{\tilde\chi_f}$ and
the universal dif\/ferential operator}

For a series $f(u) = \sum\limits_{s=0}^\infty f_su^{-s-1}$ with coef\/f\/icients in $\C$
we have an automorphism
\[
\tilde\chi_f \ :\ \glx\ \to\ \glx ,\qquad
e_{ab}^{\{s\}}\,\mapsto\,e_{ab}^{\{s\}}+ f_s\,\delta_{ab}.
\]
For a $\glx$-module $V$, denote by $V^f$ the representation
of $\glx$ on the same vector space given by the rule
\[
X|_{V^f} = (\tilde\chi_f(X))|_V
\]
for any $X\in\glx$.

\begin{lemma}
Let $\La = (\La_1,\dots,\La_N)$ be a dominant integral $\glN$-weight and
$a$ a complex number. Denote
$\La(a) = (\La_1 + a , \ldots, \La_N + a)$.
Let $V = M_\La(z)$ be an evaluation module and
$f(u)=\frac{a}{u-z}$. Then $V^f = M_{\La(a)}(z)$.
\end{lemma}


\noindent
{\bf 11.1.1.}
Let $\bs \La = (\La^{(1)},\dots,\La^{(n)})$ be a collection of
integral dominant $\glN$-weights, where $\La^{(i)}=(\La_1^{(i)},\dots,
\La_N^{(i)})$ for $i=1,\dots, n$.
Let $\bs a = (a_1,\dots,a_n)$ be a collection on complex numbers.
Introduce the new collection
$\bs\La(\bs a)=(\La^{(1)}(a_1),\dots,\La^{(n)}(a_n))$.

For complex numbers $z_1,\dots,z_n$, consider
the tensor products of evaluation modules
\begin{gather*}
\M_{\bs \La}(\bs z) = M_{\La^{(1)}}(z_1) \otimes \dots
\otimes M_{\La^{(n)}}(z_n) ,
\\
\M_{\bs \La(\bs a)}(\bs z) =
M_{\La^{(1)}(a_1)}(z_1) \otimes\dots\otimes M_{\La^{(n)}(a_n)}(z_n) .
\end{gather*}

\begin{theorem}\label{thm comparison Gaudin*}
For any matrix $K$, an $M_{\bs\La}(\bs z)$-valued function $g(u)$ belongs to
the kernel of the differential operator $\Dg_{K,\M_{\bs\La}(\bs z)}(u,\p)$ if
and only if the function $C(u)g(u)$ belongs to the kernel of the differential
operator $\Dg_{K,\M_{\bs\La(\bs a)}(\bs z)}(u,\p)$, where
\[
C(u) = \prod_{i=1}^n\, (u-z_i)^{a_i} .
\]
\end{theorem}


\noindent
{\bf 11.1.2.}
Assume that $\M_{\bs\La}(\bs z)$ is the tensor product of polynomial
$\glx$-modules and $\bs a = (a_1, \dots, a_N)$ are non-negative integers.
Then $C(u)$ is a polynomial and $\M_{\bs\La(\bs a)}(\bs z)$ is the tensor
product of polynomial $\glx$-modules. In this case
Theorem~\ref{thm comparison Gaudin*} allows us to compare the quasi-exponential
kernels of the dif\/ferential operators $\Dg_{K,\,\M_{\bs\La}(\bs z)}(u,\p)$ and
$\Dg_{K,\,\M_{\bs\La(\bs a)}(\bs z)}(u,\p)$.

\section{The kernel in the tensor product of evaluation modules\\
in the Gaudin case}
\label{The kernel of D Gaudin}

\subsection{The second main result in the Gaudin case}
Let $K\,=\,\diag\,(K_1,\dots,K_N)$ be a diagonal matrix. Let $\bs \La =
(\La^{(1)},\dots,$ $\La^{(n)})$ be a collection of integral dominant
$\glN$-weights, where $\La^{(i)}=(\La_1^{(i)},\dots, \La_N^{(i)})$ and
$\La_N^{(i)}=0$ for $i=1,\dots, n$.

Choose a weight subspace
$\M_{\bs \La}(\bs z)[\bs m]\subset M_{\bs \La}(\bs z)$ and consider
the universal dif\/ferential operator
$\Dg_{K,\,\M_{\bs \La}(\bs z)[\bs m]}(u,\p)$ associated with $K$
and $\M_{\bs\La}(\bs z)[\bs m]$.

\begin{theorem}
\label{thm second main, part 1 Gaudin}
Assume that the numbers $K_1,\dots,K_N$ are distinct.
Then for any $i=1,\dots,N$ and any nonzero vector
$v_0 \in \M_{\bs \La}(\bs z)[\bs m]$, there exists
an $\M_{\bs \La}(\bs z)[\bs m]$-valued quasi-exponential
\begin{gather}
\label{formula of a quasi-exp Gaudin}
e^{K_iu}(v_0 u^{m_i}+v_1 u^{m_i-1}+\dots+ v_{m_i})
\end{gather}
which lies in the kernel of
$\Dg_{K,\M_{\bs \La}(\bs z)[\bs m]}(u,\p)$.
Moreover, the $\C$-span of all such quasi-exponen\-tials is the kernel of
$\Dg_{K,\,\M_{\bs \La}(\bs z)[\bs m]}(u,\p)$.
\end{theorem}

\begin{proof}
On one hand, by Theorem \ref{The first main thm Gaudin}, the kernel of
$\Dg_{K,\M_{\bs\La}(\bs z)[\bs m]}(u,\p)$ is generated by quasi-exponentials.
On the other hand, by Lemma \ref{lem of Q Gaudin} and part (i) of
Theorem \ref{thm on coefficients Gaudin}, any quasi-exponential lying in the
kernel of $\Dg_{K,\,\M_{\bs \La}(\bs z)[\bs m]}(u,\p)$ must be of the form
\Ref{formula of a quasi-exp Gaudin}, where
$v_0 $ is a nonzero vector. Since
the kernel is of dimension $N \cdot \dim \M_{\bs \La}(\bs z)[\bs m]$,
there exists a quasi-exponential of the form
\Ref{formula of a quasi-exp Gaudin}
with an arbitrary nonzero $v_0$.
\end{proof}

\begin{theorem}
\label{sec local result Gaudin}
For $i=1,\dots, n$, $j=1,\dots,N$, and any nonzero vector
$v \in \M_{\bs \La}(\bs z)[\bs m]$,
there exists an $\M_{\bs \La}(\bs z)[\bs m]$-valued
solution $f(u)$ of the differential equation
$\Dg_{K,\M_{\bs \La}(\bs z)[\bs m]}(u,\p) f(u)= 0$ such that
\[
f(u) = v\,(u-z_i)^{\La^{(i)}_j + N-j} + O\big((u-z_i)^{\La^{(i)}_j + N-j+1}\big).
\]
\end{theorem}

The proof easily follows from part (ii) of
Theorem \ref{thm on coefficients Gaudin}.


\medskip

\noindent
{\bf 12.1.1.}
Assume that $K$ is the zero matrix.
Consider the subspace $\sing \M_{\bs \La}(\bs z)[\bs m] \!\subset\!
\M_{\bs \La}(\bs z)[\bs m]$ of $\glN$-singular vectors and the associated
universal dif\/ferential operator
$\Dg_{\!K{=}0,\sing\M_{\bs \La}(\bs z)[\bs m]}(u,\!\p\!).\!$

\begin{theorem}
\label{thm on space for K=0}
For any $i=1,\dots,N$ and any nonzero vector
$v_0 \in \sing\,\M_{\bs \La}(\bs z)[m_1,\dots,m_N]$, there exists a
$\sing\,\M_{\bs \La}(\bs z)[\bs m]$-valued polynomial
\begin{gather}
\label{formula for polynom Gaudin}
v_0\,u^{m_i + N-i} + v_1\,u^{m_i + N-i-1} +\cdots + v_{m_i+N-i}
\end{gather}
which lies in the kernel of
$\Dg_{K=0, \sing\M_{\bs \La}(\bs z)[\bs m]}(u,\p)$.
Moreover, the $\C$-span of all such polynomials is the kernel of
$\D_{K=0, \sing\,\M_{\bs \La}(\bs z)[\bs m]}(u,\p)$.
\end{theorem}

\begin{proof}
On one hand, by Theorem \ref{The first main thm Gaudin}, the kernel of
$\Dg_{K=0,\,\sing\M_{\bs \La}(\bs z)[\bs m]}(u,\p)$
is generated by quasi-exponentials.
On the other hand, by part (i) of
Theorem \ref{thm on coefficients Gaudin}, any quasi-exponential lying in the
kernel of $\Dg_{K=0,\,\sing\M_{\bs \La}(\bs z)[\bs m]}(u,\p)$
must a polynomial.
By Theorem \ref{thm K=0} such a polynomial has to be of the form indicated in
\Ref{formula for polynom Gaudin}, where
$v_0 $ is a nonzero vector. Since
the kernel is of dimension $N \cdot \dim\, \sing\,\M_{\bs \La}(\bs z)[\bs m]$,
there exists a polynomial of the form
\Ref{formula for polynom Gaudin}
with an arbitrary nonzero $v_0$.
\end{proof}

\subsection[Kernel of the fundamental differential operator
associated with an eigenvector]{Kernel of the fundamental dif\/ferential operator
associated\\ with an eigenvector}
\label{cor on fund spaces Gaudin}

Assume that $v \in \M_{\bs \La}(\bs z)[\bs m]$
is an eigenvector of all transfer matrices,
\[
\Tee_{K,\M_{\bs \La}(\bs z)[\bs m]}(u)\, v =
\lambda_{k,v}(u) \,v ,
\qquad
k=0,\dots, N .
\]
Then the scalar dif\/ferential operator
\[
\Dg_v(u,\p) = \sum_{k=0}^N\,(-1)^k\,\lambda_{k,v}(u)\,\p^{N-k}
\]
will be called {\it the fundamental differential operator}
associated with the eigenvector $v$.

Theorems \ref{thm second main, part 1 Gaudin}, \ref{sec local result Gaudin}
and~\ref{thm on space for K=0} give us information on the kernel of the
fundamental dif\/ferential operator.

\begin{corollary} \ \
\begin{enumerate}\itemsep=0pt
\item[(i)]
Assume that $K$ is diagonal with distinct diagonal entries.
Then the kernel of $\Dg_v(u,\p)$ is generated by quasi-exponentials
$e^{K_1u}p_1(u),\dots, e^{K_Nu}p_N(u)$, where for every $i$
the polynomial $p_i(u)\in \C[u]$ is of degree $m_i$.

\item[(ii)]
Assume that $K$ is the zero matrix and the eigenvector $v$ belongs to
$\sing\,\M_{\bs \La}(\bs z)[\bs m]$. Then the kernel of $\Dg_v(u,\p)$
is generated by suitable polynomials $p_1(u),\dots, p_N(u)$ of degree
$m_1+N-1,\dots, m_N$, respectively.
\end{enumerate}
\end{corollary}

\begin{corollary}
For $i=1,\dots, n$, $j=1,\dots,N$, there exists a solution $f(u)$
of the differential equation $\Dg_v(u,\p)\,f(u)\,=\,0$ such that
\[
f(u) = (u-z_i)^{\La^{(i)}_j + N-j} + O\big((u-z_i)^{\La^{(i)}_j + N-j+1}\big).
\]
\end{corollary}

If $v$ is a Bethe eigenvector, then these two corollaries are proved in
\cite{MV1} and \cite{MV3}.

\subsection*{Acknowledgments}

E.~Mukhin was supported in part by NSF grant DMS-0601005.
V.~Tarasov was supported in part by RFFI grant 05-01-00922.
A.~Varchenko was supported in part by NSF grant DMS-0555327.

\pdfbookmark[1]{References}{ref}

\LastPageEnding
\end{document}